\documentclass[12pt]{amsart}
\font\emailfont=cmtt10

\usepackage{amsmath,amsthm,amsfonts,amscd,flafter,epsf}

\newcommand\commentable[1]{#1}

\newcommand{\rk}{\mathrm{rk}}

\newtheorem{theorem}{Theorem}[section]
\newtheorem{prop}[theorem]{Proposition}
\newtheorem{cor}[theorem]{Corollary}

\newtheorem{lemma}[theorem]{Lemma}

\newtheorem{defn}[theorem]{Definition}

\def\endproof{\relax\ifmmode\expandafter\endproofmath\else
  \unskip\nobreak\hfil\penalty50\hskip.75em\hbox{}\nobreak\hfil\bull
  {\parfillskip=0pt \finalhyphendemerits=0 \bigbreak}\fi}
\def\endproofmath$${\eqno\bull$$\bigbreak}
\def\bull{\vbox{\hrule\hbox{\vrule\kern3pt\vbox{\kern6pt}\kern3pt\vrule}\hrule}}

\newcommand{\R}{\mathbb{R}}

\newcommand{\Z}{\mathbb{Z}}

\newcommand{\OneHalf}{\frac{1}{2}}

\newcommand{\cm}{\cdot}

\newcommand{\ModSWfour}{\mathcal{M}}
\newcommand{\ModFlow}{\ModSWfour}

\newcommand\sgn{\mathrm{sgn}}

\newcommand\Hom{\mathrm{Hom}}

\newcommand\abuts\Rightarrow
\newcommand\Sym{\mathrm{Sym}}

\newcommand\HFpRed{\HFp_{\red}}
\newcommand\HFpRedEv{\HFp_{\red, \ev}}
\newcommand\HFpRedOdd{\HFp_{\red,\odd}}

\newcommand{\ev}{\mathrm{ev}}
\newcommand{\odd}{\mathrm{odd}}

\newcommand\relspinc{\underline{\spinc}}

\newcommand\x{\mathbf x}

\newcommand\y{\mathbf y}

\newcommand\ModSphere{\ModFlow\left({\mathbb S}\longrightarrow 
\Sym^{g-1}(\Sigma_{1})\times \Sym^2(\Sigma_{2})\right)}
\newcommand\ModSpheres\ModSphere

\newcommand\CFa{\widehat{CF}}

\newcommand{\red}{\mathrm{red}}

\newcommand\HFp{\HFb}

\newcommand\HFa{\widehat{HF}}
\newcommand\HFb{HF^+}

\newcommand\UnparModSp{\widehat \ModSp}
\newcommand\UnparModFlow\UnparModSp
\newcommand\Mod\ModSp

\newcommand{\cald}{{\mathcal D}}

\newcommand{\spinc}{\mathfrak s}

\newcommand\ModMaps{\mathcal M}
\newcommand\ModSp\ModMaps

\newcommand\alphas{\mbox{\boldmath$\alpha$}}

\newcommand\betas{\mbox{\boldmath$\beta$}}

\newcommand\Dual{\mathcal D}
\newcommand\Duality\Dual

\newcommand\essarc{\boldmath{E}}
\newcommand\essarcP{\essarc_+}
\newcommand\essarcM{\essarc_-}

\newcommand\spincrel\relspinc

\newcommand\CFK{CFK}
\newcommand\HFK{HFK}

\newcommand\CFKa{\widehat\CFK}

\newcommand\HFKa{\widehat\HFK}

\newcommand\BasePt{w}
\newcommand\FiltPt{z}

\commentable{  

\title[{Knot Floer homology, genus bounds, and mutation}] 
{Knot Floer homology, genus bounds, and mutation}

\author[Peter Ozsv{\'a}th]{Peter Ozsv\'ath}
\address{Department of
Mathematics, Columbia University, New York 10025 \newline
\indent{\emailfont{petero@math.columbia.edu}}}
\thanks{PSO was partially supported by NSF grant number DMS-0234311}

\author[Zolt{\'a}n Szab{\'o}]{Zolt{\'a}n Szab{\'o}} 
\address{Department of
Mathematics, Princeton University, New Jersey 08544 \newline
\indent{\emailfont{szabo@math.princeton.edu}}}}
\thanks{ZSz was partially supported by NSF grant numbers DMS-0107792
and a Packard Fellowship.}

\begin{document}

\begin{abstract}
In an earlier paper, we introduced a collection of graded Abelian
groups $\HFKa(Y,K)$ associated to knots in a three-manifold.  The aim
of the present paper is to investigate these groups for several
specific families of knots, including the Kinoshita-Terasaka knots and
their ``Conway mutants''.  These results show that $\HFKa$ contains
more information than the Alexander polynomial and the signature of
these knots; and they also illustrate the fact that $\HFKa$ detects
mutation. We also calculate $\HFKa$ for certain pretzel knots, and
knots with small crossing number ($n\leq 9$).  Our calculations 
give obstructions to certain Seifert fibered surgeries on the
knots considered here.
\end{abstract}

\maketitle
\section{Introduction}

In~\cite{Knots}, we defined an invariant for knots $K\subset
S^3$, which 
take the form
of a graded Abelian group $\HFKa(K,i)$ for each 
integer $i$. The main results of~\cite{AltKnots} give explicit descriptions of
some of the input required for determining $\HFKa$ in terms of the
combinatorics of a generic planar projection of $K$. As an
application, it is shown that $\HFKa$ for an alternating knot is
explicitly determined by the Alexander polynomial and the signature of
the knot (compare also~\cite{Rasmussen}).  The aim of the present
article is to apply and extend techniques from~\cite{AltKnots} to determine certain
knot homology groups of some more complicated types of knots.
Indeed, to underscore the relative strength of $\HFKa$ over the
Alexander polynomial, we focus mainly on certain knots with trivial
Alexander polynomial (and hence vanishing signature). 

These calculations have the following consequences.  Of course, they
show that $\HFKa$ is stronger than the Alexander polynomial; but more
interestingly, they also show that, unlike many other knot invariants,
$\HFKa$ is sensitive to Conway mutation. These computations
further underline an interesting relationship between the 
knot Floer homology and the Seifert genus $g(K)$
of the knot $K$. Specifically,
recall that in Theorem~\ref{Knots:thm:Adjunction} of~\cite{Knots}, we proved an adjunction
inequality, stating that 
if $\deg \HFKa(K)$ denotes the largest integer $d$ for which
$\HFKa(K,d)\neq 0$, then
\begin{equation}
\label{eq:Adjunction}
\deg\HFKa(K) \leq g(K).
\end{equation}
Indeed, we also conjectured,  based on 
the analogy with Seiberg-Witten theory and a theorem of 
Kronheimer and Mrowka~\cite{KMThurston},
that 
$$\deg\HFKa(K) =g(K)$$
for every knot in $S^3$. 
Calculations from this paper
can be taken as further evidence supporting this conjecture.

Finally, the calculations provide obstructions to realizing Seifert
fibered spaces as certain surgeries on $S^3$ along many of the knots studied
here.

We emphasize that in general, calculating $\HFKa$ is not a purely
combinatorial matter. The generators of this complex can be described
combinatorially, and indeed in~\cite{AltKnots}, we indentified them
with Kauffman states (c.f.~\cite{Kauffman}), but the differentials count
pseudo-holomorphic disks in a symmetric product. However, there are
some additional combinatorial aspects of this chain complex described
below (see Section~\ref{sec:Heegaards}), including a multi-filtration
on the chain complex, which facilitate our calculations. As a further 
illustration of these techniques, we also calculate the knot Floer homology groups
for all knots with at most nine crossings.

We now give a description of the knots we study
and state the results of our calculations. 
 
\subsection{Kinoshita-Terasaka and Conway knots}

In~\cite{KinoshitaTerasaka}, Kinoshita and Terasaka construct a family
of knots $KT_{r,n}$, indexed by integers $|r|\neq 1$ and $n$, with
trivial Alexander polynomial.  These knots are obtained by modifying a
picture of the $(r+1,-r,r,-r-1)$ (four-stranded) pretzel links, and
introducing $2n$ twists.  There are some redundancies in these knots.
When $r\in \{0,1,-1,-2\}$ or $n=0$, this construction gives the
unknot.  Also, there is a symmetry identifying $KT_{r,n}=KT_{-r-1,n}$,
which can be realized by turning the knot inside out. Finally, the
reflection of $KT_{r,n}$ is the knot $KT_{r,-n}$. Now, recall that
the knot Floer homology groups transform in a controlled manner under
reflection: i.e. if ${\overline K}$ denotes the reflection of $K$,
then for each $i,d\in\Z$, $$\HFKa_d(K,i)\cong \HFKa^{-d}({\overline
K},-i)$$ (where here the left-hand-side denotes knot Floer homology in
dimension $d$, while the right-hand-side denotes knot Floer {\em
co}-homology in dimension $-d$) and also
\begin{equation}
\label{eq:Symmetry}
\HFKa_d(K,i) \cong \HFKa_{d-2i}(K,-i)
\end{equation}
(c.f. Equations~\eqref{Knots:eq:Reflection} and~\eqref{Knots:eq:KnotConjInv}
respectively of~\cite{Knots}),
so there is no loss of
generality in assuming $r>1$ and $n>0$.  We have illustrated the case where
$r=3$ and $n=2$ in Figure~\ref{fig:KT}. It is possible
to eliminate one crossing from the diagram for $KT_{r,n}$, but the
new diagram is somewhat more cumbersome to draw.

We calculate the topmost non-trivial knot Floer homology group for
$KT_{r,n}$ in Section~\ref{sec:CalcKT}, arriving
at the following result:

\begin{figure}
\mbox{\vbox{\epsfbox{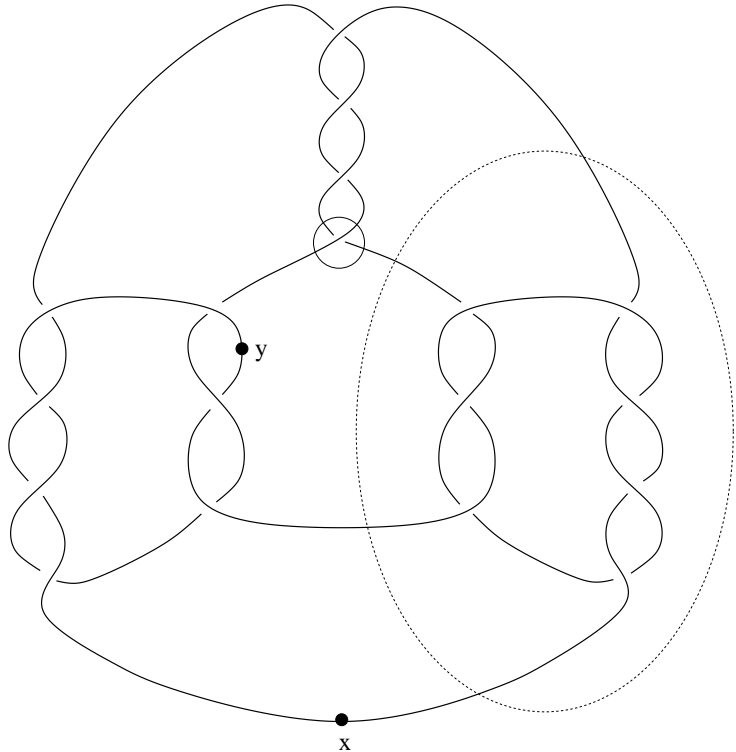}}}
\caption{\label{fig:KT}
{\bf{Kinoshita-Terasaka knot with $r=3$ and $n=2$.}}  When the circled
crossing is changed, we obtain $KT_{3,1}$, while if it is resolved, we
obtain a $(r+1,-r,r,-r-1)$-pretzel link. The Conway knot is obtained
as a mutation around the sphere indicated here with a large, dotted
ellipse (indeed, it is the mutation induced by $180^\circ$ rotation
about the axis perpendicular to the plane of the knot projection).
The relevance of the indicated points $x$ and $y$ will become apparent
in Sections~\ref{sec:CalcKT} and \ref{sec:CalcConway} respectively.}
\end{figure}

\begin{theorem}
\label{thm:KTcalc}
Consider the Kinoshita-Terasaka knot $KT_{r,n}$ with $n>0$ and
$r>1$.  This knot has $\HFKa(KT_{r,n},s)=0$ for all $s>r$, and
$$
\HFKa(KT_{r,n},r)\cong
\Z^n_{(r)} \oplus \Z^n_{(r+1)},
$$
where here (and indeed throughout this paper)
the subscript $(r)$ indicates that the
corresponding summand is supported in dimension $r$.
\end{theorem}

Note that in~\cite{GabaiKT}, Gabai exhibits a Seifert surface for
$KT_{r,n}$ with genus $r$, and proves that it is genus-minimizing, using
the theory of foliations.
It is interesting to note that Theorem~\ref{thm:KTcalc}, together
with Inequality~\eqref{eq:Adjunction}, gives an
alternate proof that this Seifert surface is genus-minimizing.
Some new applications will be described later (c.f. Subsection~\ref{subsec:NoSeif}).

Theorem~\ref{thm:KTcalc} is based on the results of~\cite{AltKnots},
where we give combinatorial descriptions of some of the data for
calculating $\HFKa$. In fact, in Section~\ref{sec:Heegaards} we
explain how some of this data can be simplified, and amplify it with a
multi-filtration on the chain complex of Kauffman states. Using these
techniques, we obtain some additional calculations, as well.

Let $C_{r,n}$ denote the Conway knot, which is obtained from $K_{r,n}$
by mutation. This knot is obtained using the same construction as
$K_{r,n}$, only using a four-stranded pretzel link of type
$(r+1,-r,-r-1,r)$ instead of $(r+1,-r,r,-r-1)$. Alternatively, it can
be thought of as obtained from $K_{r,n}$ by a mutation using the
sphere pictured in Figure~\ref{fig:KT}, c.f.~\cite{Lickorish}.  These
knots also have trivial Alexander polynomial, and indeed,
they satisfy the same symmetries as $K_{r,n}$. Note that these knots, too,
admit a projection with one fewer crossing. In the case where
$r=2$ and $n=1$, an eleven-crossing projection is pictured
in
Figure~\ref{fig:DomainConway}. We prove the following in
Section~\ref{sec:CalcConway}:

\begin{theorem}
\label{thm:Conway}
Let $C_{r,n}$ denote the Conway mutant of $KT_{r,n}$ with $n>0$ and $r>1$.
This knot has $\HFKa(C_{r,n},s)=0$ if $s>2r-1$, and
$$\HFKa(C_{r,n},2r-1)\cong 
\Z_{(2r-1)}^n\oplus \Z_{(2r)}^n. 
$$
\end{theorem}

It is easy to construct Seifert surfaces $F$ for $C_{r,n}$
with genus $g(F)=2r-1$, compare~\cite{GabaiLinks}, see also
Section~\ref{sec:CalcConway}. 

Since $KT_{r,n}$ and $C_{r,n}$ differ by a Conway mutation, and their
groups $\HFKa$ are manifestly different, we see that, unlike the
Alexander, Jones, HOMFLY, and Kauffman polynomials, the invariant
$\HFKa$ is sensitive to mutation. It is interesting to compare this
with Khovanov's invariants, c.f. ~\cite{KhovMut}, \cite{Khovanov}.

\subsection{Pretzel knots}

The techniques described here also lend themselves quickly
to a calculation for pretzel knots $P(p,q,r)$,
where $p$, $q$, and $r$ are odd integers.
When $p$, $q$, and $r$ all have the
same sign, these knots are alternating, and hence their Floer homology
has been determined in~\cite{AltKnots}.  Thus, by reflecting the knot
if necessary and relabeling, we are left with the case where $q<0$ and
$p,r>0$.

\begin{theorem}
\label{thm:PretzCalc}
Consider the knot $K=P(p,q,r)$ where $p=2a+1$, $q=-(2b+1)$, $r=2c+1$,
with $a,b,c\geq 0$. Then, if $b\geq \min(a,c)$, we have that 
$$\HFKa(K,1)=\Z^{ab+bc+b-ac}_{(1)}.$$
If $b<\min(a,c)$,
we have that
$$\HFKa(K,1)=\Z^{b(b+1)}_{(1)}\oplus \Z^{(b-a)(b-c)}_{(2)}.$$
\end{theorem}

This family contains infinitely many knots with trivial Alexander
polynomial: the Alexander polynomial is trivial precisely when
$pq+qr+pr+1=0$ (e.g. let $(p,q,r)=(-3,5,7)$). It follows at once from
the above theorem that for all non-trivial pretzel knots in the above
family, $\HFKa$ is also non-trivial.

\subsection{Knots with few crossings}
\label{subsec:SmallKnots}

Although the techniques from Section~\ref{sec:Heegaards} are not
sufficient to calculate $\HFKa$ in general, they can be employed
successfully in the study of relatively small knots, as measured by
the number of double-points. In fact, in Section~\ref{sec:SmallKnots}
we calculate $\HFKa$ for all knots with nine or fewer crossings,
except for two particular knots whose $\HFKa$ has been calculated
in~\cite{Knots} and \cite{NoteLens} (the knot $9_{42}$ and $8_{19}$ --
the $(3,4)$ torus knot). The Floer homology of the remaining
non-alternating knots (with less than ten crossings) behaves like the 
Floer homology of
alternating knots, c.f. Theorem~\ref{thm:SmallKnots} below.

\subsection{Surgeries on knots}
\label{subsec:NoSeif}

In another direction, the calculations of this paper can be used to
give information on three-manifolds obtained as  surgeries on
knots, following results of~\cite{SurgSeif} (see
also~\cite{SomePlumbs}).  To explain, recall that the formal sum of
Euler characteristics of $\HFKa$ gives the symmetrized Alexander
polynomial $\Delta_K(T)$:
\begin{equation}
\label{eq:Euler}
\sum_{i} \chi \left(\HFKa(K,i)\right) \cm T^i =\Delta_K(T)
\end{equation}
(c.f. Section~\ref{Knots:sec:Links} of~\cite{Knots}).
It is an immediate corollary of this that
\begin{equation}
\label{eq:AlexDefect}
\deg\Delta_K\leq \deg\HFKa(K).
\end{equation}

It is a result of~\cite{SurgSeif} (see
esp. Corollary~\ref{SurgSeif:cor:AlexHFKDefect} of~\cite{SurgSeif})
that if $K$ is a knot for which $\deg\HFKa(K)>1$ and
Inequality~\eqref{eq:AlexDefect} is strict, then $K$ 
does not admit certain Seifert fibered surgeries. 
Specifically,  we have the following:

\begin{cor}
For any integer $q\neq 0$, $1/q$ surgery on $S^3$ along $KT_{n,r}$ or
$C_{n,r}$ (with $n>0$ and $r>1$) is never a Seifert fibered space.
\end{cor}

For the case of pretzel knots $P(p,q,r)$ with $p$, $q$, and $r$ all
odd, Corollary~\ref{SurgSeif:cor:AlexHFKDefect} of~\cite{SurgSeif} no
longer applies, since $\deg\HFKa(K)=1$. And indeed, there are cases of
such pretzel knots with Seifert fibered surgeries. However, a careful
look at the proof of that corollary, and a closer look at $\HFKa$
gives the following corollary
(c.f. Proposition~\ref{prop:NoSeifSurg}). Note that this corollary
covers all non-trivial three-stranded pretzel knots with trivial
Alexander polynomial (compare with~\cite{GodaTer} and~\cite{Mattman}):

\begin{cor}
\label{cor:NoSeifSurgPretz}
Let $P(p,q,r)$ be a non-trivial pretzel knot with $p$, $q$, and $r$
odd. When $p\geq r>0$ and $q<-1$ with $|q|< \min(p-2,r)$, no integral
surgery along $P(p,q,r)$ is a Seifert fibered space.
\end{cor}

\noindent{\bf{Further remarks.}} Additional calculations of 
knot Floer homology groups can be found in~\cite{RasmussenThesis}
and~\cite{Eftekhary}.  The authors wish to thank Eaman Eftekhary,
Cameron Gordon,
Mikhail Khovanov, Rob Kirby, Paul Melvin, and Jacob Rasmussen for
interesting conversations.

\section{Calculational tools}
\label{sec:Heegaards}

Let $K\subset S^3$ be a knot.  In~\cite{Knots}, we introduced the knot
Floer complex $\CFKa(K)=\bigoplus_{s\in\Z}\CFKa(K,s)$ which is
associated to a Heegaard diagram for a knot, and whose homology groups
are knot invariants, see
also~\cite{RasmussenThesis}. In~\cite{AltKnots}, we gave a description
of the generators of the chain complex $\CFKa$ in terms of
combinatorics of a generic projection for a knot (together with some
extra data). We recall the constructions in
Subsection~\ref{subsec:Simplify}, and show that in some cases, the
number of generators can be cut down, to make the calculations
simpler.  In Subsection~\ref{subsec:Domains}, we give a combinatorial
description of the domain of a homotopy class $\phi\in\pi_2(\x,\y)$
connecting a pair of states. Here, the condition that $\cald(\phi)\geq
0$ from~\cite{HolDisk} (a necessary condition for $\y$ to appear with 
non-zero multiplicity in the expression for $\partial \x$)
is formulated in terms of a multi-filtration on
the set of states.

In Subsection~\ref{subsec:Trees}, we recall the correspondence between
states and maximal subtrees (see also~\cite{Kauffman}), which removes
much of the redundancy which is inherent in the description of a
state. 

In Subsection~\ref{subsec:Skein}, we turn to certain properties of
the ``skein exact sequence'' from~\cite{Knots}.

\subsection{Simplifying Heegaard diagrams}
\label{subsec:Simplify}

Choose an orientation for $K$ and a generic projection of $K$ to the
plane.  The projection gives a planar graph $G$ where the vertices of
$G$ correspond to the double-points of the projection of $K$, and the
edges inherit an orientation from $K$.  Choose a distinguished edge $\epsilon_0$
for this planar graph. We call a projection
with this additional data a {\em decorated knot projection}.

There are four distinct quadrants (bounded by edges) emanating from
each vertex, each of which is a corner of the closure of
some region of
$S^2-G$. We distinguish the two of these regions which contain the
distinguished edge on their boundary, denoting them ${\mathbf A}$ and ${\mathbf B}$.

\begin{defn}
A {\em state} (c.f.~\cite{Kauffman}) is an
assignment which associates to each vertex of $G$ one of the four
in-coming quadrants, so that:
\begin{itemize}
\item the quadrants associated to distinct vertices are subsets of distinct
regions in $S^2-G$
\item none of the quadrants associated to vertices 
is a corner of the distinguished regions
${\mathbf A}$ or ${\mathbf B}$.
\end{itemize}
\end{defn}

It is easy to see that a state sets up a one-to-one correspondence
between vertices of $G$ and the regions of $S^2-G-{\mathbf A}-{\mathbf B}$. 

\begin{figure}
\mbox{\vbox{\epsfbox{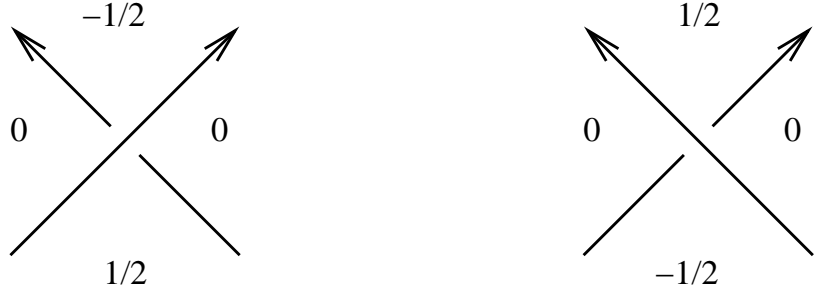}}}
\caption{\label{fig:FiltLevel}
{\bf{Local filtration level contributions.}} We have
illustrated the local contributions for the
filtration level of a state for both kinds of
crossings.}
\end{figure}

\begin{figure}
\mbox{\vbox{\epsfbox{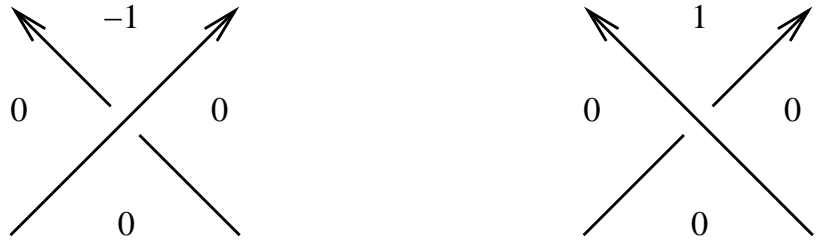}}}
\caption{\label{fig:AbsGrading}
{\bf{Local grading contributions.}} We have
illustrated the local contribution for the absolute grading
associated for a state.}
\end{figure}

\begin{defn}
\label{def:FiltGrad}
The {\em filtration level} of a state is the integer obtained by 
adding up the local contributions at each
quadrant, which are determined by the crossing types of the knots, as
indicated in Figure~\ref{fig:FiltLevel}. 
The {\em absolute grading} of a state is
determined by adding up another local contribution, pictured in
Figure~\ref{fig:AbsGrading}.
\end{defn}

We have the following result from~\cite{AltKnots}:

\begin{theorem}
\label{thm:States}
Fix a generic knot projection for a knot $K\subset S^3$.  There is a
one-to-one correspondence with generators of $\CFKa(K,i)$ and states
$x$ whose associated filtration level is $i$. Under this
correspondence, the absolute grading of a state in
the above sense coincides with the
absolute degree of the corresponding generator of $\CFKa(K,i)$.
\end{theorem}

\begin{defn}
\label{def:EssentialInt}
Fix a decorated knot projection. An {\em essential interval} is a sequence of consecutive
edges
with the following properties:
\begin{itemize}
\item the distinguished edge $\epsilon_0$ appears in the sequence, so
we can write
$$\essarc= \bigcup_{i=-\ell}^m\epsilon_i,$$
where here $\ell, m \geq 0$, and
$\epsilon_{i+1}$ is the successor of $\epsilon_i$ for all $i=-\ell,...,m-1$,
\item the immersed arcs 
\begin{eqnarray*}
\essarcP=\bigcup_{i=1}^m\epsilon_i &{\text{and}}&
\essarcM=\bigcup_{i=-\ell}^{-1}\epsilon_i
\end{eqnarray*}
are pairwise disjoint
\item as we traverse the arc $\essarcP$ according to the orientation of $K$,
i.e. starting at the vertex $\epsilon_0\cap \epsilon_1$, and then
passing through $\epsilon_1,...,\epsilon_m$ in order, 
all of the crossings we encounter for the first time have the same type
(i.e. they are all either over- or under-crossings);
similarly, as we traverse the arc $\essarcM$
backwards, i.e. starting at $\epsilon_0\cap\epsilon_{-1}$,
and proceeding through till $\epsilon_{-\ell}$, all the crossings we encounter
the first time have the same type (which might be different different from the
crossing type encountered along $\essarcP$).
\end{itemize}
\end{defn}

\begin{defn}
\label{def:Essential}
Fix a decorated knot projection, and also an essential interval
$\essarc=\bigcup_{i=-\ell}^m\epsilon_i$.  Each edge $\epsilon_i$
inherits an orientation from $K$, and we write its endpoints as
$\partial \epsilon_i=v_{i+1}-v_i$, so that $\{v_{-\ell},...,v_{m+1}\}$
are the vertices in the order they appear along $\essarc$.  An
$\essarc$-essential state is a state $x$ with the following
properties:
\begin{itemize}
\item for $i=1,...,m$, if $v_i\not\in \{v_1,v_2,...,v_{i-1}\}$, 
then the corner
containing $x(v_i)$ contains the edge $\epsilon_i$ on its boundary
\item for $i=-\ell+1,...,0$, 
if $v_i\not\in\{v_{0},v_{-1},...,v_{i+1}\}$, then the corner containing
$x(v_i)$ contains the edge $\epsilon_{i-1}$ on its boundary.
\end{itemize}
\end{defn}

\begin{prop}
\label{prop:Simplify}
Fix a knot projection $G$ for a knot $K\subset S^3$, and let $\essarc$
be an essential arc. Then, there is a system of generators for $\CFKa$
consisting of only $\essarc$-essential states (in the sense of
Definition~\ref{def:Essential}), with gradings and filtration levels
as given in Theorem~\ref{thm:States}.
\end{prop}

\begin{proof}
Recall that the Heegaard surface constructed in~\cite{AltKnots} is
obtained as a boundary regular neighborhood of the knot
projection. For each vertex $v$, we have a $\beta$-curve denoted
$\beta_v$, and at the distinguished $\epsilon_0$, we choose a meridian
$\mu$ for the knot which is supported near $\epsilon_0$. Then, on
either side of that meridian in the Heegaard surface, we choose a pair
of basepoints $\BasePt$ and $\FiltPt$ for the definition of
$\CFKa(K)$.

Suppose for simplicity that all the $v_i$
are distinct (i.e. that $\essarcM\cup\essarcP$ is an embedded arc in $G$).
To simplify the Heegaard diagram as above, we move the two basepoints
further away from the meridian, so that we can handleslide the
$\beta$-curves belonging to $v_{1}$ and $v_{0}$ across $\mu$, to get
new curves $\beta_{1}'$ and $\beta_{0}'$. 
We continue handlesliding in
this manner -- $\beta_{v_{k+1}}$ across $\beta_{v_k}'$ and
$\beta_{v_{-k}}$ across $\beta_{v_{-k+1}}'$ -- until we have a new
sequence of $\beta$-curves $\beta_{v_0}',...,\beta_{v_k}'$ which now
meet always at most two $\alpha$-curves (rather than four). The
hypothesis on the crossing types was used to ensure that the reference
points could always be moved out of the handlesliding region (without
crossing any of the attaching circles). 
This procedure is illustrated
in Figure~\ref{fig:SimplifyDiag}.
It is not difficult to modify the above procedure when $\essarcM$ and $\essarcP$ have
self-intersections.

\begin{figure}
\mbox{\vbox{\epsfbox{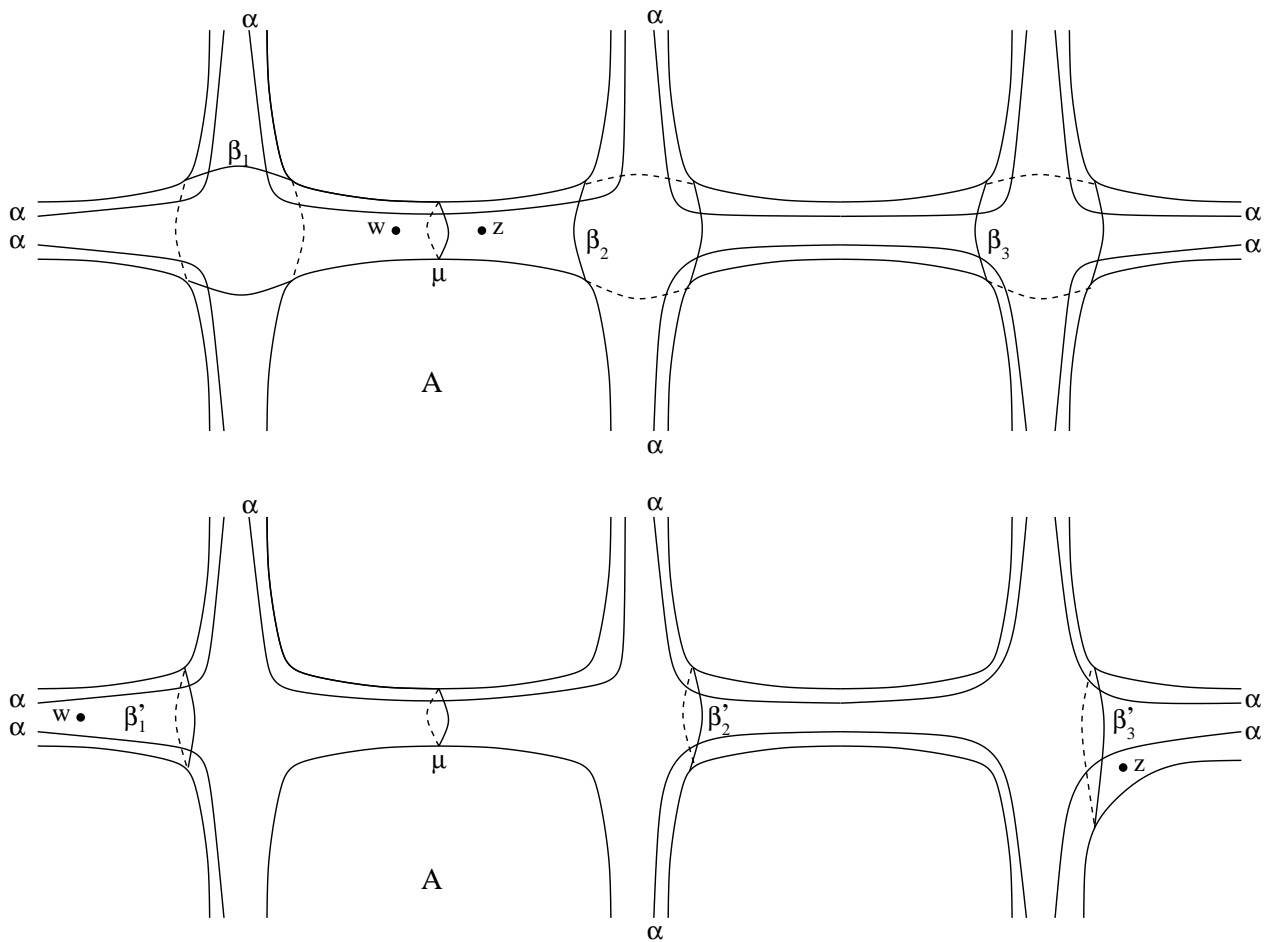}}}
\caption{\label{fig:SimplifyDiag}
{\bf{Diagram simplification.}} We have illustrated the proof of
Proposition~\ref{prop:Simplify}.  At the top, we have illustrated the
part of the Heegaard diagram coming from~\cite{AltKnots} (near where
we have two under-crossings followed by one over-crossing), choosing
our reference point between the under- and the over-crossing.  We have
dropped all the subscripts to the $\alpha$-curve, but we remind the
reader that there is one region in each of the compact components of
$\R^2-G$ (the non-compact region here is denoted by ${\mathbf A}$).
At the bottom, we have
illustrated the corresponding ``simplification''.}
\end{figure}
\end{proof}

Of course, the notion of essential arc, and the above proposition depends on the choice of 
essential interval. However, for some decorated knot
projections, there is always a unique maximal
essential interval. Indeed, this is the case for all the projections we consider in this paper, and
hence, with this understood, we call a state an essential state if it is $\essarc$-essential
for this maximal essential interval $\essarc$.

\subsection{The combinatorics of domains, and the multi-filtration}
\label{subsec:Domains}

Of course, the calculation of $\HFKa(K)$ requires an explicit
understanding of the differential in the complex $\CFKa(K)$, which in
turn involves a pseudo-holomorphic curve count. More precisely, given a
pair of states $x$ and $y$ (whose filtration level coincides, and
whose absolute gradings differ by one), letting
$\x$ and $\y$ denote the corresponding generators of $\CFKa(K)$,
there is a unique homotopy
class of Whitney disk $\phi\in\pi_2(\x,\y)$ with
$n_\FiltPt(\phi)=n_\BasePt(\phi)=0$. The $\y$ component
of $\partial \x$  is the signed count of points $\# \ModFlow(\phi)/\R$
in the moduli space of pseudo-holomorphic representatives of
$\phi$ (modulo translation). At present, this count does not have a
direct combinatorial description.  However, there are combinatorial
conditions on $x$ and $y$ which ensure that it vanishes.

Let $G$ be a graph for the knot projection of $K$ with $N$ edges. We
order the edges $\{\epsilon_i\}_{i=0}^{N-1}$ of $G$, so that
$\epsilon_0$ is the marked edge, and the others appear in the order in
which they are encountered by moving along $K$ (with its specified
orientation).  Let $v_i$ denote the vertex at the intersection of
$\epsilon_i$ with $\epsilon_{i+1}$ (note that each vertex in $G$
appears as $v_i$ for two different values of $i$, once as an
overcrossing, and once as an undercrossing).  We can associate to each
state $x$ a multi-filtration-level
$$M_x\in\Hom(\{\epsilon_i\}_{i=0}^{N-1},\Z\oplus \Z),$$ as follows: $$
M_x(\epsilon_i)=
\left\{
\begin{array}{ll}
(0,0) & {\text{if $i=0$}} \\
M_x(\epsilon_{i-1})+(0,1) & {\text{if $v_i$ is an overcrossing and $x(v_i)$
is to the right of $\epsilon_{i}\cup \epsilon_{i-1}$}} \\
M_x(\epsilon_{i-1})-(0,1) & {\text{if $v_i$ is an overcrossing and $x(v_i)$
is to the left of $\epsilon_{i}\cup \epsilon_{i-1}$}} \\
M_x(\epsilon_{i-1})+(1,0) & {\text{if $v_i$ is an undercrossing and $x(v_i)$
is to the left of $\epsilon_{i}\cup \epsilon_{i-1}$}} \\
M_x(\epsilon_{i-1})-(1,0) & {\text{if $v_i$ is an undercrossing and $x(v_i)$
is to the right of $\epsilon_{i}\cup \epsilon_{i-1}$}} 
\end{array}\right.$$

\begin{defn}
\label{def:PO}
Fix a decorated knot projection. We define a partial ordering on the set
of states, as follows: if $x$ and $y$ are two different states, then
$x>y$ if for all edges $\epsilon$, $M_x(\epsilon)-M_y(\epsilon)$ is a
pair of non-negative integers.
\end{defn}

\begin{figure}
\mbox{\vbox{\epsfbox{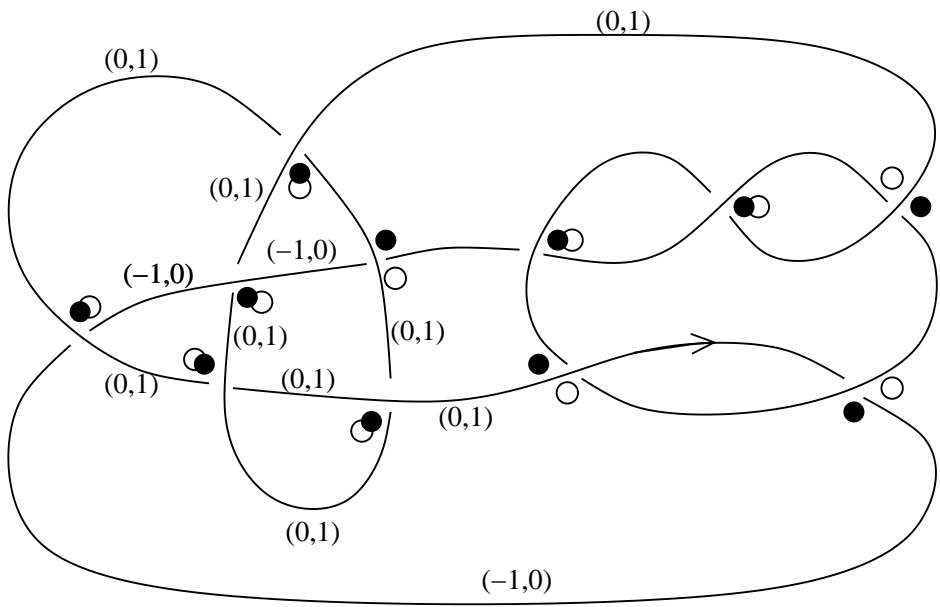}}}
\caption{\label{fig:DomainConway}
{\bf{A domain.}}
Let $x$ denote the state indicated by the dark circles and $y$ denote
the state indicated by the hollow ones, in the Conway knot as
pictured above. Then, $\OneHalf(M_x-M_y)$ is represented as above: specifically,
near each edge $\epsilon$, we have written $\OneHalf(M_x-M_y)(\epsilon)$,
unless the latter quantity vanishes. Since sometimes this is negative,
it follows that $x\not> y$, and hence, although it is easy to see
that $\deg(x)=\deg(y)+1$, $y$ does not appear in the expansion for
$\partial x$.}
\end{figure}

An example of two states $x$ and $y$ for which $x\not>y$ and $y\not>x$
is illustrated in Figure~\ref{fig:DomainConway}.

\begin{prop}
\label{prop:Domains}
Suppose that $x$ and $y$ represent the same filtration level, and indeed 
suppose that $y$ appears in $\partial x$ with non-zero multiplicity,
then $x>y$.
\end{prop}

\begin{proof}
Let $(\Sigma,\alphas,\betas)$ denote the Heegaard diagram for $S^3$
used for Theorem~\ref{thm:States}. Consider an edge $\epsilon_i$ in
the knot projection which does not meet the one of the two
distinguished regions in the knot diagram. This edge, then, gives rise
to a cylinder in the Heegaard diagram, which is divided into two
squares by the $\alpha$-arcs. We place one reference point $t_i$ in
the ``top'' part of the diagram, and another one $b_i$ in the
``bottom'' part. It is easy to see (after a straightforward
case-by-case analysis of $x(v_i)$ and $y(v_i)$) that if
$\phi\in\pi_2(\x,\y)$ is the homotopy class with
$n_\FiltPt(\phi)=n_\BasePt(\phi)=0$, then
$$\OneHalf\Big(M_x(\epsilon_i)-M_y(\epsilon_i)\Big)
=(n_{t_i}(\phi), n_{b_i}(\phi)).$$ 
The result now follows from the basic fact that if $\phi$
has a pseudo-holomorphic representative (for suitably small
perturbations of the holomorphic condition), then all these
multiplicities must be non-negative (c.f. Lemma~\ref{HolDisk:lemma:NonNegativity} 
of~\cite{HolDisk}).
\end{proof}

Note that the above argument also applies when the Heegaard diagram is
``simplified'' as in Proposition~\ref{prop:Simplify}, only in that
case, one uses the mult-filtration only over those edges which are not
in the essential interval.

\subsection{Trees}
\label{subsec:Trees}

States admit a rather economical description in graph-theoretic terms
(see also~\cite{Kauffman}).  The regions in the complement of the
planar projection can be colored black and white in a chessboard
pattern, by the rule that any two regions which share an edge have
opposite color. There is then an associated ``black graph'', whose
vertices correspond to the regions colored black, and whose edges
correspond to vertices in $G$, which connect the opposite black
regions. We let ${\mathbf A}$ (resp. ${\mathbf B}$) denote the black (resp. white)
region whose boundary contains the distinguished edge $\epsilon_0$.

In these terms, states are in one-to-one correspondence with the
maximal subtrees of the black graph.  Given a state $x$, we associate
to it the union of vertices of $G$, thought of now as edges in the
black graph, to which $x$ associates a black quadrant. This gives the
maximal black subtree associated to the state $x$.

Conversely, given a black subtree $T$, we can orient the edges so that
the ``root'' is the distinguished black region ${\mathbf A}$, and all edges
point away from this root. We construct the black part of the
corresponding vertex assignment, as follows. Let $v$ be a
vertex of $G$ which corresponds to some edge of $T$. With respect to
the induced orientation on $T$, this oriented edge of $T$ points to a
uniquely determined endpoint $r\in T$, which in turn corresponds to
one of the two black quadrants (in $S^2-G$) which meet at $v$. We
let $x(v)$, then, be the quadrant corresponding to $r$. To determine
the rest of the vertex assignment, we first consider the dual white
graph $T^*$, obtained from the white graph by deleting all the edges
corresponding to the vertices appearing in the black subtree $T$. Note
that $T^*$ is actually a tree, and repeat the above procedure, now for
the white quadrants.

Clearly, one can reformulate the results of Theorem~\ref{thm:States}
in terms of these graphs. Each edge of both the black and white graphs
inherits a label among the numbers $\{-1, 0 , +1\}$: an edge in the
white (resp. black) graph is labelled with $0$ if the two white
(resp. black) quadrants meeting at the corresponding vertex both are
have grading contribution $0$ (c.f. Figure~\ref{fig:AbsGrading}), and it
is labelled with $\pm 1$ if one of the  two white (resp. black)
graph is labelled with $\pm 1$. We call edges labelled with $0$ {\em neutral} edges, and view them as unoriented, while the non-neutral ones are
oriented (so as to point away from the vertex corresponding
to the quadrant labelled with $0$ in Figure~\ref{fig:AbsGrading}).

Given a black tree $T$, let $T^*$ be its dual white tree, and orient
both as before (so that the edges of $T$ and $T^*$ point away from
${\mathbf A}$ and ${\mathbf B}$ respectively). Then, twice the
filtration level corresponding to $T$ is obtained as follows. We sum
over all edges $e$ in $T$ the label of the edge $e$ in $G$ times $+1$
if its orientation as an edge of the black graph agrees with the the
induced orientation coming from $T$, and $-1$ if the orientations
disagree (note here that this last sign is irrelevant for neutral
edges); and then add to that the corresponding sum for all edges of
the dual white graph $T^*$.

Similarly, the grading level corresponding to $T$ is obtained by summing
over all edges $e$ in $T$ the label of the edge $e$ in $G$, provided
that its orientation agrees with the induced orientation coming from $T$
(and zero otherwise), and once again adding the corresponding 
sum for $T^*$.

\subsection{Exact Sequences}
\label{subsec:Skein}

In a different direction, we derived in~\cite{Knots} a ``skein exact
sequence'', which we recall here.  

Suppose that $K$ is an oriented knot in $Y$, and we have a disk $D$
which meets $K$ in two algebraically cancelling points. If we perform
$-1$ surgery on $\gamma=\partial D$, we obtain a new knot $K_-$ in
$S^3$, which is obtained from $K$ by introducing a full twist in a
tubular neighborhood along $\gamma$. Performing $+1$ surgery on
$\gamma$ gives us another knot $K_+$ with a full twist introduced in
the other direction.  There is a third link $K_0$ which is obtained by
resolving the knot (so as to miss $D$ entirely). Recall~\cite{Knots}
that the link invariant in this case is, by definition, the knot
invariant for the knot in $S^1\times S^2$ obtained by performing a
zero-surgery along $\gamma$.

In~\cite{Knots}, we established skein exact sequences for each integer $i$
$$
\begin{CD}
...@>>>\HFKa(K,i)@>{f_1}>>\HFKa(K_0,i)@>{f_2}>>\HFKa(K_+,i) @>{f_3}>> ...
\end{CD}
$$
and 
$$
\begin{CD}
...@>>>\HFKa(K_-,i)@>{g_1}>>\HFKa(K_0,i)@>{g_2}>>\HFKa(K,i) @>{g_3}>> ...,
\end{CD}
$$
where the maps $f_i$ and $g_i$ are induced by two-handle additions.

It will be useful to us to have the following compatibility result about
$g_2$ and $f_1$:

\begin{lemma}
\label{lemma:CompositeLemma}
For the above two exact sequences, the composite
$$g_2\circ f_1\colon \HFKa(K,i) \longrightarrow \HFKa(K,i)$$
is trivial.
\end{lemma}

\begin{figure}
\mbox{\vbox{\epsfbox{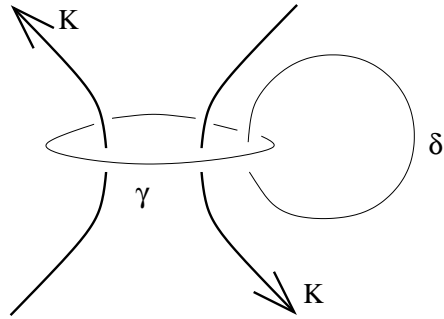}}}
\caption{\label{fig:CompositeMap}
{\bf{Composite map.}}  We have pictured here two strands of $K$.
Performing zero-surgery on $\gamma$ realizes the map $f_1$, while zero-surgery
on $\delta$ realizes $g_2$.}
\end{figure}

\begin{proof}
The map $f_1$ is induced by cobordism formed by zero-surgery on
$\gamma$, while $g_2$ is induced by the cobordism formed by
zero-surgery on another unknot $\delta$ which links $\gamma$ once (and
does not link the knot $K$).  As in~\cite{HolDiskFour}, we obtain the
same map if we switch the order in which we perform the two two-handle
additions: first we perform zero-surgery on $\delta$, and then
zero-surgery on $\gamma$ (see Figure~\ref{fig:CompositeMap}).
However, in this latter ordering, we factor through the three-manifold
$Y\#(S^2\times S^1)$ (and the knot is contained entirely in the $Y$
summand. It is easy to adapt the K\"unneth principle for connected
sums in this context (see especially
Proposition~\ref{HolDiskTwo:prop:ConnSum} of~\cite{HolDiskTwo}) to see
that 
\begin{equation}
\label{eq:Kunneth}
\HFKa(Y\#(S^2\times S^1),K,i)\cong\HFKa(Y,K,i)\otimes
H^1(S^2\times S^1).
\end{equation}
Of course, $H^1(S^2\times S^1)\cong \Z\oplus
\Z$, but our reason for writing the answer in this form is that we
have now an action of $[\delta]\in H_1(S^2\times S^1)$ on
$\HFKa(Y\#(S^2\times S^1),i)$ (induced from the corresponding action
on $\HFa(Y\#(S^2\times S^1))$). The above isomorphism is compatible
with this action of $[\delta]$ (where $[\delta]\in H_1(S^2\times S^1)$
acts on the right-hand-side through its natural action on $H^1(S^2\times S^1)$, c.f. 
Section~\ref{HolDisk:subsec:DefAct} of~\cite{HolDisk}).
Now, $f_1$ maps to
the kernel of the $[\delta]$-action, since the curve representing
$\delta$ is null-homologous in the four-manifold obtained by attaching
a two-handle along $\delta$ to $Y-K$; while $g_2$ is trivial on the
image of the $[\delta]$-action because, once again, $\delta$ is
null-homologous in the corresponding four-manifold
(compare~\cite{HolDiskFour}).  But in the description of the
$[\delta]$-action from Equation~\eqref{eq:Kunneth}, 
it is clear that the image of $[\delta]$ coincides
with the kernel of $[\delta]$, completing the proof.
\end{proof}

\section{Calculations for the Kinoshita-Terasaka knots}
\label{sec:CalcKT}

Consider the Kinoshita-Terasaka knots $KT_{r,n}$ with $n\neq 0$.  All
the knots $KT_{r,n}$ have trivial Alexander polynomial, but the knots
themselves are non-trivial when $r>1$ and $n\neq 0$. Indeed,
Gabai exhibits a Seifert surface for $KT_{r,n}$ of genus $r$ and
proves (c.f.~\cite{GabaiKT}) that this Seifert surface has minimal
genus.

\begin{figure}
\mbox{\vbox{\epsfbox{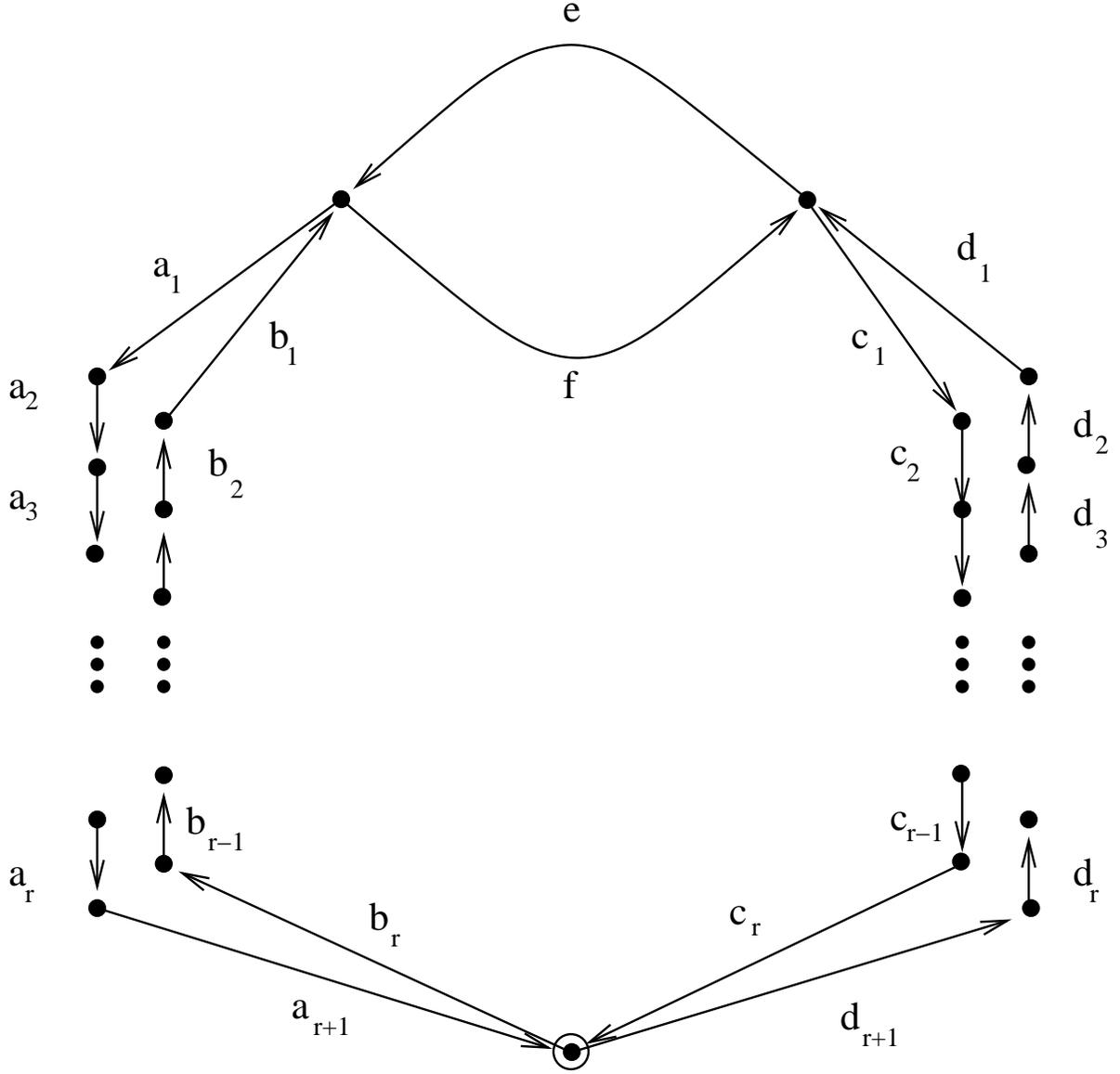}}}
\caption{\label{fig:KTdiag}
{\bf{Black graph for the Kinoshita-Terasaka knot.}}  The following
edges are labeled with $+1$: $e$, $f$ $\{b_i\}_{i=1}^{r}$,
$\{d_i\}_{i=1}^{r+1}$; the following are labeled with $-1$:
$\{a_i\}_{i=1}^{r+1}$, $\{c_i\}_{i=1}^{r}$.  The distinguished black
vertex is circled. The white regions correspond to the complement of
this graph in $S^2$, and the distinguished white region is the
non-compact one (i.e. the one bounded by the chain of $\{a_i\}$,
$\{d_i\}$, and $e$).  Note that all the edges in the dual white graph
are neutral (since none of the edges in the black graph are).}
\end{figure}

We shall focus first on the case where $n=1$. We distinguish the edge
connecting the base of the $-r-1$ and $r+1$ strands opposite to
where the ($n=1$) twisting takes place -- this is indicated by the point $x$ pictured
in Figure~\ref{fig:KT}.  We have illustrated the ``black graph'' of
$KT_{r,1}$, in Figure~\ref{fig:KTdiag}. For the black graph, there are
two edges, labeled $e$ and $f$, and four chains $\{a_i\}_{i=1}^{r+1}$,
$\{b_i\}_{i=1}^{r}$, $\{c_i\}_{i=1}^{r}$ and
$\{d_i\}_{i=1}^{r+1}$. The edges $e$, $f$, $b_i$, and $d_i$ are
labeled with $+1$, while those of type $a_i$ and $c_i$ are labeled
with $-1$. The distinguished black region corresponds to the vertex
where $a_{r+1}$, $b_r$, $c_r$, and $d_{r+1}$ meet.  The vertices of
the white graph are, of course, the regions in the complement of this
planar graph. The region bounded by the chain of $b_i$ and $c_i$ and
the edge $f$ corresponds to the distinguished white region. In fact,
since none of the black edges is neutral, all of the white ones are.

We calculate $\HFKa(KT_{r,n},r)$ in the case where $n=1$ using
Theorem~\ref{thm:States}, together with the multi-filtration
(Proposition~\ref{prop:Domains}).  The case where $n$ is arbitrary
will follow from properties of the skein exact sequence
c.f. Lemma~\ref{lemma:CompositeLemma} (though one could alternately
give a more direct argument using the multi-filtrations).

\begin{figure}
\mbox{\vbox{\epsfbox{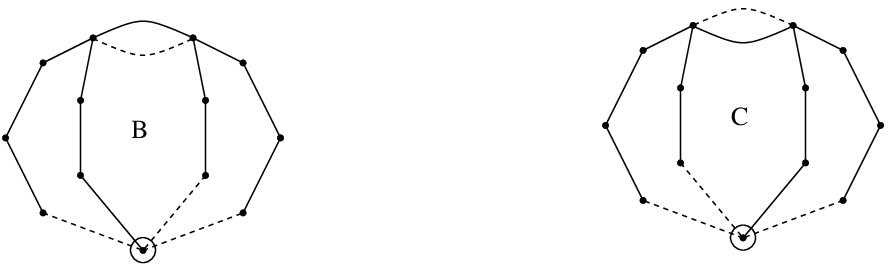}}}
\caption{\label{fig:KTTrees}
{\bf{Generators for $K_{n,r}$ with $n=1$, $r=2$.}}}
\end{figure}

\vskip.2cm
\noindent{\bf{Proof of Theorem~\ref{thm:KTcalc} when $n=1$.}}
It is easy to see that there are two trees $B$ and $C$ which represent
filtration level $-r$: the tree $B$ does not contain $a_{r+1}$, $f$,
$c_r$, or $d_{r+1}$, while the tree $C$ does not contain $a_{r+1}$,
$e$, $b_r$, or $d_{r+1}$. These trees are illustrated in Figure~\ref{fig:KTTrees}.

\begin{figure}
\mbox{\vbox{\epsfbox{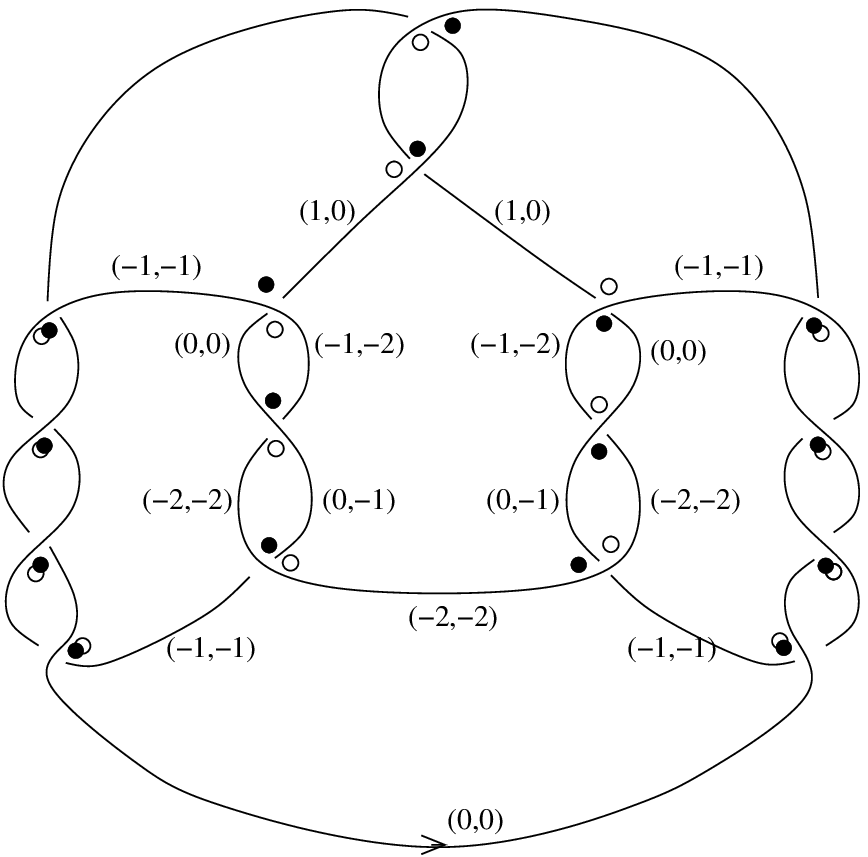}}}
\caption{\label{fig:KTDomain}
{\bf{Comparison of generators $B$ and $C$.}}
We have labeled some of the multiplicities for the multi-filtration
comparing these two generators of 
$\HFKa(KT_{n,r},-r)$, when $n=1$ and $r=2$. Here, $B$ is
denoted by the black dot, and $C$ by the hollow one.
The arrow appears on
the distinguished edge of the knot projection (at the vertex denoted
by $x$ in Figure~\ref{fig:KT}).}
\end{figure}

Moreover, the grading of $B$ is given by
$1-r$, while the grading of $C$ is given by $-r$.  It is
straightforward to verify that $B\not>C$ for the multifiltration 
(see Figure~\ref{fig:KTDomain} for an illustration)
and
hence, according to Proposition~\ref{prop:Domains}, $\partial\equiv
0$. These calculations show that
$$\HFKa(KT_{r,1},s)\cong 
\left\{\begin{array}{ll}
0 & {\text{if $s<-r$}} \\
\Z_{(-r)}\oplus \Z_{(1-r)} & {\text{if $s=-r$.}}
\end{array}
\right.$$
Note that this is equivalent to the statement of the theorem (with
$n=1$), in view of the symmetry 
of the knot Floer homology groups 
(c.f. Equation~\eqref{eq:Symmetry}).

\vskip.2cm
\noindent{\bf{Calculation for an $(r+1,-r,r,-r-1)$ pretzel link}.}
Note that there is a skein exact sequences of the form
\begin{equation}
\label{eq:SkeinKT}
\begin{CD}
...
@>>>\HFKa(KT_{r,n-1},i) @>{f^{n}_1}>>\HFKa(X(r),i) @>{f^{n}_2}>> \HFKa(KT_{r,n},i) @>>>...,
\end{CD}
\end{equation}
where $X(r)$ is an oriented $(r+1,-r,r,-r-1)$ pretzel link which is
independent of $n$. (Strictly speaking, the notation for $f^n_1$ and
$f^n_2$ ought to include the level $i$, but we suppress this
for readibility.) The maps $f^n_1$ and $f^n_2$
both decrease absolute grading by $\OneHalf$. In fact, according to
Lemma~\ref{lemma:CompositeLemma},
\begin{equation}
\label{eq:CompositeLemma}
f^{n}_2\circ f^{n+1}_1\equiv 0.
\end{equation}
Moreover, $KT_{r,0}$ is the unknot. Specializing to the case where
$i=r$, and $n=1$ it follows at once that
\begin{equation}
\label{eq:PretzelLink}
\HFKa(X(r),r)\cong \Z_{(r-\OneHalf)}\oplus \Z_{(r+\OneHalf)},
\end{equation}
and the map $f^1_2$ is an isomorphism. 

\vskip.2cm
\noindent{\bf{Proof of Theorem~\ref{thm:KTcalc} for arbitrary $n$.}}
To establish the theorem for all $n\geq 1$ we prove inductively both
the theorem, and also the statement that the map $f_2^n$ is
injective. The basic case was established above.
For the inductive step, if $f_2^{n+1}$ is not injective,
then $f_1^{n+1}$ would have to be
non-trivial, but this contradicts the injectivity of $f_2^{n}$
(which holds by the inductive hypothesis), together with
Lemma~\ref{lemma:CompositeLemma}, in the form of
Equation~\eqref{eq:CompositeLemma}. 
\qed

\subsection{Additional remarks}

Consider the case of $K_{r,n}$ with $n=1$.  By moving the marked edge,
we obtain various chain complexes representing
$\HFKa(K_{r,1})$. Typically, subtrees which represent filtration level
$r$ vary as we move the marked edge. However, since the white graph
consists of neutral edges only, if we choose our marked edge so
that the distinguished black region $X$ remains unchanged -- there are
four possible choices -- then the maps from maximal subtrees to
filtration levels and degrees are unchanged. However, the map from
maximal subtrees to states, of course, varies, and more
interestingly, the induced partial ordering on 
subtrees can change, too. For example,
if we mark the edge of $B$ opposite to the edge containing $x$
(c.f. Figure~\ref{fig:KT}), then it is easy to see that the two
generators $B$ and $C$ described in the proof of
Theorem~\ref{thm:KTcalc} are still the two representatives for
filtration level $r$, and their dimensions are $r+1$ and $r$. However,
with this choice of marked edge, it is now the case that $B>C$.

\section{Calculations for the Conway knots}
\label{sec:CalcConway}

We consider the Conway mutants $C_{r,n}$ of the Kinoshita-Terasaka
knots.  The calculation of $\HFKa(C_{r,n},2r-1)$ proceeds similarly to
the calculations from Section~\ref{sec:CalcKT}.  Note that the knot
$C_{2,1}$ can be given the eleven-crossing presentation pictured in
Figure~\ref{fig:DomainConway}.  In this case, if we place the
reference point where the arrow is indicated, and use the
simplification of the Heegaard diagram described in
Section~\ref{sec:Heegaards}, then it is straightforward to see that
in filtration level $3$ 
there are only two inessential states, the two states $x$ and $y$
pictured in the figure, and they have absolute grading $4$ and $3$
respectively. However, $\partial x=0$, since $x\not>y$, as illustrated
in the figure, verifying Theorem~\ref{thm:Conway} for $r=2$, $n=1$.

As before, we begin by restricting to the case where $n=1$.  Rather
than drawing the black graph in this case, we indicate the necessary
modifications to Figure~\ref{fig:KTdiag}. (Note that it is not hard to
find projections with fewer essential states than the ones we describe
-- indeed, in the desired filtration level, we can arrange for there
to be only two essential states, as in the case with $r=2$ considered
above. However the diagrams we describe presently have the advantage
that they are easier to describe in words.)  The black graph of
$C_{r,n}$ looks just like that for $KT_{r,n}$, except that now there
are $r+1$ edges of type $c_i$, and only $r$ edges of type
$d_i$. Moreover, the edges labeled with $+1$ are $e$, $f$,
$\{b_i\}_{i=1}^r$ and $\{c_i\}_{i=1}^{r+1}$, and the edges labeled
with $-1$ are $\{a_i\}_{i=1}^{r+1}$ and $\{d_i\}_{i=1}^{r+1}$. In
particular all edges in the dual white graph remain neutral.  We
choose our marked edge to contain the point $y$ in Figure~\ref{fig:KT}
(after mutating). Correspondingly, now, the distinguished black edge
at the vertex between $b_1$ and $b_2$, and the distinguished white
region is bounded by the chains $\{b_i\}_{i=1}^r$,
$\{c_i\}_{i=1}^{r+1}$ and the edge $f$.

\begin{figure}
\mbox{\vbox{\epsfbox{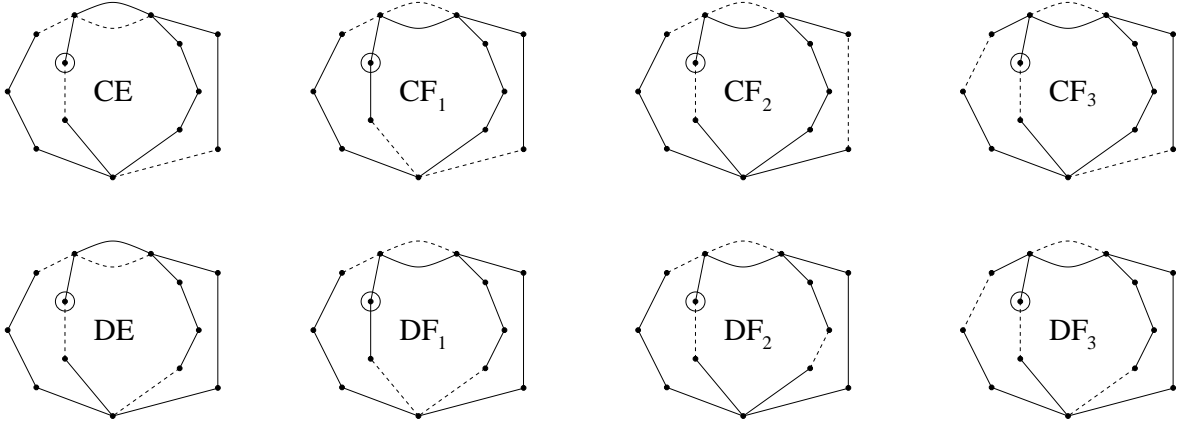}}}
\caption{\label{fig:TreesConway}
{\bf{Generators for $C_{n,r}$ with $n=1$, $r=3$.}}}
\end{figure}

\begin{lemma}
\label{lemma:EnumGensConway}
For the diagram of $C_{r,1}$ described here with $r>2$, there are $8$
Kauffman states representing filtration level $2r-1$. Of these, $3$
are in dimension $r+2$, which we denote $CE$, $CF_1$, $CF_2$, and
$CF_3$, and $3$ are in dimension $r+1$, which we denote $DE$, $DF_1$
$DF_2$, and $DF_3$.  (Note that these states are explicitly identified in the
proof below, which also explains their notation. See also
Figure~\ref{fig:TreesConway}.) In the case where $r=2$, we have
only $6$ essential states: the two states $CF_1$ and $DF_1$ are missing
from this diagram.
\end{lemma}

\begin{proof}
Notice that there are a total of $4r+4$ edges in the black graph, and
in fact a maximal spanning tree must contain exactly $4r$ edges.
It is clear that a maximal tree which contains neither $e$ nor $f$
has filtration level at most $r$. 

Thus, all maximal subtrees with filtration level $2r-1$ contain either
$e$ or $f$. Now, if $T$ is a maximal subtree which contains one of $E$
or $F$, then clearly, it must contain at least one of the chains
$A=\{a_i\}$, $B=\{b_i\}$, $C=\{c_i\}$ or $D=\{d_i\}$. Indeed, it
cannot contain more than one. Accordingly, we say that a tree is of
type $AE$ if it contains the chain $A$ and the vertex $e$. We analyze the
eight cases separately.

Trees of type $AE$ which represent filtration level $2r-1$, we claim,
cannot contain the edge $b_1$. But a tree which contains $A$, and does
not contain $b_1$ corresponds to an inessential vertex assignment --
and indeed, in the case where $r>2$, such trees never represent
filtration level $2r-1$.  The same remarks rule out trees of type
$AF$, and similar remarks rule out trees of type $BE$, and $BF$.

It is easy to see that there is only one tree of type $CE$ which
represents filtration level $2r-1$, and it is the one which does not
contain $a_1$, $f$, $b_2$, and $d_{r}$. Similarly, there is only one tree
of type $DE$, and it does not contain $a_1$, $f$, $b_2$, and $c_{r+1}$.

Assume for the moment that $r>2$.  There are three trees of type $CF$
representing filtration level $2r-1$: one which does not contain $e$,
$a_1$, $b_3$, and $d_{r}$, which we denote $CF_1$, one which does not
contain $e$, $a_1$, $b_2$, and $d_2$, which we denote $CF_2$, and one
which does not contain $e$, $a_2$, $b_2$, and $d_{r}$. Similarly,
there are three trees of type $DF$ representing filtration level
$2r-1$, which we denote $\{DF_i\}_{i=1}^3$, where $DF_i$ is gotten
from $CF_i$ by deleting $e$ and $c_{r+1}$ and adding $f$ and $d_r$.
The case where $r=2$ works similarly, except that the states $CF_1$
and $DF_1$ do not exist.
\end{proof}

\begin{lemma}
\label{lemma:OrderGensConway}
Consider the essential generators listed in 
Lemma~\ref{lemma:EnumGensConway}
for $C_{r,1}$ with $r>2$.
These have the following ordering properties for the multi-filtration:
\begin{equation}
\begin{array}{llll}
CF_1 >& CF_2  >& CE  > &CF_3 \label{eq:Ordering} \\
DF_1 >& DF_2  >&DE  > &DF_3,
\end{array}
\end{equation}
and
\begin{equation}
\begin{array}{lll}
CF_3\not> DE, & CE\not> DF_2, & CF_2\not> DF_2 \label{eq:NoOrdering} \\
\end{array}
\end{equation}
Moreover, fixing fix $i\neq 2$, and letting $\phi\in\pi_2(CF_i,DF_i)$
be the homotopy class with $n_\BasePt(\phi)=0$, we have that
\begin{equation}
\label{eq:BoundaryMap}
\#\left(\frac{\ModFlow(\phi)}{\R}\right)=1.
\end{equation}
The same holds for the
corresponding homotopy class in $\pi_2(CE,DE)$.
When $r=2$, the same remarks hold, excluding states $CF_1$ and $DF_1$.
\end{lemma}

\begin{proof}
Verifying the order properties is straightforward.  Some of the work
in verifying Relation~\eqref{eq:NoOrdering} is shortened, given
Relation~\eqref{eq:Ordering}, and the observation that the difference
in the multi-filtering between $CF_i$ and $DF_i$ for $i\neq 2$ is
supported in a single edge; a similar remark holds for $CE$ and $DE$.
For the statement about homotopy classes, observe that all of these
homotopy classes are represented by quadrilaterals, and hence
$\#\ModFlow(\phi)/\R$ can be calculated by one-variable complex
analysis (compare Section~\ref{HolDiskTwo:sec:Examples} of~\cite{HolDiskTwo}).
\end{proof}

\vskip.2cm
\noindent{\bf{Proof of Theorem~\ref{thm:Conway}.}}
Again, we start with the case where $n=1$.
We use the complex described in Lemma~\ref{lemma:EnumGensConway}.

According to Lemma~\ref{lemma:OrderGensConway}
(c.f. Relations~\ref{eq:Ordering}) together with the basic 
property of the multi-filtration (Proposition~\ref{prop:Domains}),
the chain complex admits a subcomplex generated by $CF_3$ and $DF_3$.
Indeed, the homology of this complex is trivial, in view of 
Equation~\eqref{eq:BoundaryMap}. Thus, $\HFKa(C_{r,1},2r-1)$
is calculated as the homology of the induced quotient complex.
Another application of this principle allows us to cancel also the generators
$CE$ and $DE$. The leftover complex, generated by $CF_1$, $CF_2$,
$DF_1$, and $DF_2$, now admits a quotient complex which is
generated by $CF_1$ and $DF_1$ and hence, according to
Equation~\eqref{eq:BoundaryMap}, has trivial homology. (Note that
this step is skipped when $r=2$.)

Thus, $\HFKa(C_{r,1},2r-1)$ is
calculated by the homology of the
remaining complex generated by $CF_2$ and $DF_2$.
Since $CF_2\not>DF_2$, it follows that the homology is
$$\Z_{(2r-1)}\oplus \Z_{(2r)},$$
verifying the calculation of $\HFKa(C_{r,n},2r-1)$
when $n=1$.

Note that for the pictured knot projection, are two states representing
filtration level $2r+1$ (and none representing higher filtration levels).
However, these two states are quickly seen to cancel (they, too,
are connected by a quadrilateral).

To go from $n=1$ to arbitrary $n$, observe that 
$C_{r,0}$ is an unknot. Thus, we have
a skein
exact sequence relating the various Conway knots, 
corresponding to Equation~\eqref{eq:SkeinCon}:
\begin{equation}
\label{eq:SkeinCon}
\begin{CD}
...
@>>>\HFKa(C_{r,n-1},i) @>{f^{n}_1}>>\HFKa(Y(r),i) 
@>{f^{n}_2}>> \HFKa(C_{r,n},i) @>>>...,
\end{CD}
\end{equation}
for any integer $i$,
where now $Y(r)$ is an oriented $(r+1,-r,-r-1,r)$ pretzel link,
rather than the pretzel link $(r+1,-r,r,-r-1)$
belonging to the Kinoshita-Terasaka
knots considered earlier. (Of course, $Y(r)$ is a mutant of $X(r)$,
and both have trivial Alexander polynomial.)

With these remarks in place,
the induction used to verify the theorem runs exactly
as it did in the case of $KT_{r,n}$. As a consequence,
we also obtain the following formula for $Y(r)$:
$$
\HFKa(Y(r),2r-1)\cong \Z_{(2r-\OneHalf)}\oplus \Z_{(2r+\OneHalf)}.
$$

\qed

\subsection{A Seifert surface}

Note that the Conway knot $C_{r,n}$ has a genus $2r-1$ Seifert
surface, obtained by a straightforward
modification of the picture for the case where $r=2$, $n=1$
described in~\cite{GabaiLinks}. 

\begin{figure}
\mbox{\vbox{\epsfbox{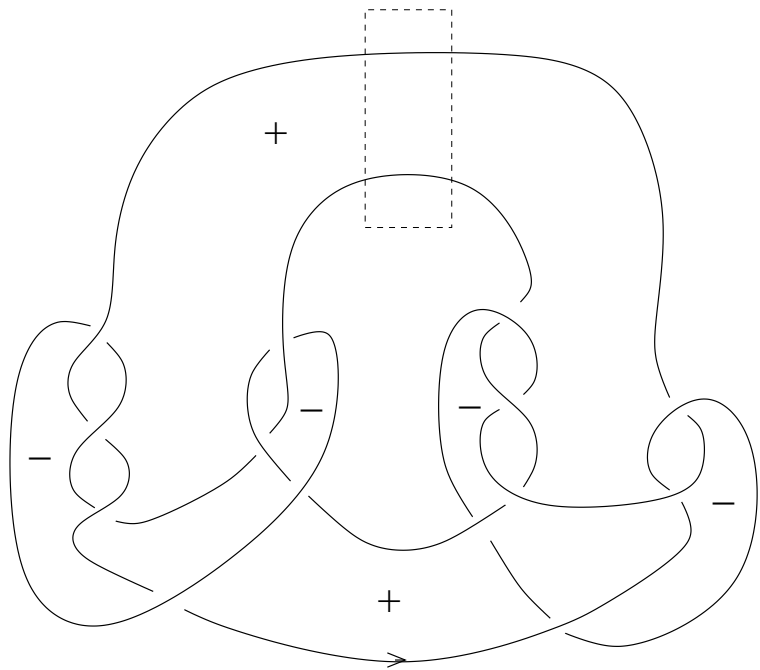}}}
\caption{\label{fig:SeifertSurface}
{\bf{A Seifert surface.}} We have illustrated here a genus $2r-2$ Seifert
surface for the $(r+1,-r,-r-1,r)$ pretzel link, when $r=2$. A Seifert
surface for the Conway knot is obtained as a Murasugi sum 
with the cylinder (with $n$ full twists) at the indicated (dashed)
rectangle.}
\end{figure}

We describe first the Seifert surface for the $(r+1,-r,-r-1,r)$
pretzel link.  We ``pull down'' the overcrossing which connects the
first two tassles and the undercrossing which connects the second
two. And then, we consider the black regions for the knot
projection. Those in turn we label with signs, with the rule that no
two regions which meet at a vertex have the same sign. The Seifert
surface for this pretzel link is obtained by connecting these black
regions by half-twists at each vertex. It is easy to see that the
surface $F$ obtained in this manner is orientable, and has
$\chi(F)=4-4r$. We have illustrated this data in
Figure~\ref{fig:SeifertSurface} for the case where $r=2$.

We can plumb this with a cylinder with $n$ full twists in it (i.e. forming a
Murasugi sum), to obtain a Seifert surface for the Conway knot whose
genus is $2r-1$.

\section{Calculations for the pretzel knots}
\label{sec:PretzCalc}

\begin{figure}
\mbox{\vbox{\epsfbox{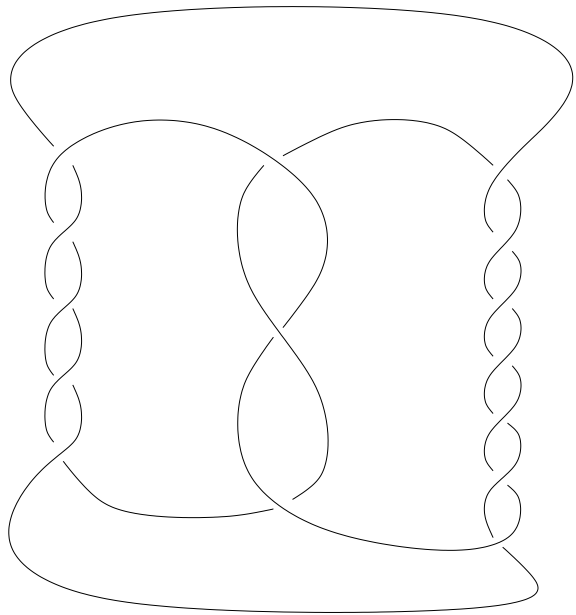}}}
\caption{\label{fig:Pretzel}
{\bf{The pretzel knot $P(5,-3,7)$.}}  This 
knot has  trivial Alexander polynomial.}
\end{figure}

In this section, we consider the family of pretzel knots $P(p,q,r)$,
where $p$, $q$ and $r$ are odd. We follow the usual conventions from
knot theory here (c.f.~\cite{Lickorish}) for the direction of the
twisting (which, unfortunately, seems to be opposite from the
convention used in~\cite{AbsGraded}), compare
Figure~\ref{fig:Pretzel}.

There are some relations amongst the pretzel knots. For example,
it is easy to see that $P(p,q,r)=P(q,r,p)$, and that
$P(-1,1,r)$ is the unknot, for any $r$. 

Recall that
$$\Delta_{P(p,q,r)}(T)=
\frac{1}{4}\Big(\left(pq+qr+pr\right)(T-2+T^{-1})
+\left(T+2+T^{-1}\right)\Big);$$ thus, there are infinitely many
examples with trivial Alexander polynomial.  When $p$, $q$, and $r$
are all positive, then the signature of $P(p,q,r)$ is given by
$$\sigma(P(p,q,r))=2.$$

When $m$ is an even integer, $T_{2,m}$ denote the torus link, oriented
so that the two strands -- which we can think of as supported inside a
solid torus -- are oriented so that the algebraic intersection of
$T_{2,m}$ with a disk transverse to the solid torus is zero. In this case,
$$\Delta_{T_{2,m}}=\frac{m}{2}(T^{-1/2}-T^{1/2}).$$ 
Moreover, the signature
of $T_{2,m}$ is $\pm 1$, depending on the sign of $m$:
$$\sigma(T_{2,m})=\sgn(m).$$

Clearly,
if we resolve one of the intersection points corresponding to the
first strand in $P(p,q,r)$, we obtain the torus link $T_{2,q+r}$.  
Thus, the skein long exact sequence of~\cite{Knots} in this case
gives 
\begin{equation}
\label{eq:SkeinPretz}
\begin{CD}
...@>>>\HFKa(P(p,q,r))@>{F}>>\HFKa(T_{2,q+r})@>{G}>> \HFKa(P(p-2,q,r))@>{H}>>...,
\end{CD}
\end{equation}
where here $F$ and $G$ preserve filtration levels, and both drop absolute 
grading by $\OneHalf$. The map $H$ also preserves filtration level,
and it preserves the parity of absolute grading.

When $p$, $q$, and $r$ all have the
same sign, then the usual projection of
$P(p,q,r)$ is alternating. In this case, 
the knot Floer homology
is  determined by~\cite{AltKnots}. Specifically, 
in~\cite{AltKnots}, it is shown that if $L$ is a
non-split, oriented, alternating link with signature
$\sigma=\sigma(L)$, and Alexander-Conway polynomial $\Delta_L(T)$,
then if we write 
$$(T^{-1/2}-T^{1/2})^{n-1}\cm \Delta_L(T)= a_0 +\sum_{s>0} a_s (T^s+T^{-s}),$$
then 
\begin{equation}
\label{eq:AltLinks}
\HFKa(L,s)\cong \Z^{|a_s|}_{(s+\frac{\sigma}{2})}.
\end{equation}

Thus, by reflecting $P(p,q,r)$ knot if necessary,
we are left with the case where $q<0$ and
$p,r>0$.  Note that when $m\neq 0$, $T_{2,m}$ is a non-split, alternating
link.

\begin{prop}
\label{prop:EasyPretzel}
Consider the pretzel knot $K=P(2a+1,-(2b+1),2c+1)$ with $a,b,c\geq 0$.
Then, if $b\geq \min(a,c)$, we have that 
$$\HFKa(K,1)=\Z^{ab+bc+b-ac}_{(1)}.$$
\end{prop}

\begin{proof}
For fixed $b\geq 0$, we prove the result by induction on both $a$ and $c$.

Consider the base case where $a=c=0$.
Using the skein exact sequence in the form of
Equation~\eqref{eq:SkeinPretz} with $p=r=1$, and the relation that
$P(-1,q,1)$ is the unknot, we see at once that
$$\HFKa_*(P(1,-(2b+1),1),1)\cong \HFKa_{*+\OneHalf}(T_{2,-2b}).$$
Moreover, since $T_{2,-2b}$ is alternating,
Theorem~\ref{AltKnots:thm:FloerHomology} of~\cite{AltKnots} applies,
and hence, in this case Equation~\eqref{eq:AltLinks} specializes to give
$$\HFKa(T_{2,-2b},1)\cong \Z^{b}_{(\OneHalf)}.$$

For the inductive step on $a$, suppose we know the result for
$P(2a+1,-2b-1,2c+1)$, and suppose that $b\geq a+1$ and $b\geq c$. The
condition that $b\geq c$ ensures that $T_{2,-2(b-c)+1}$ still has
signature $-1$, and hence $\HFKa(T_{2,-2(b-c)+1},1)$ is supported in
dimension $1/2$; by the inductive hypothesis,
$\HFKa(P(2a+1,-2b-1,2c+1),1)$ is supported in dimension $1$. Thus, the
skein exact sequence forces $\HFKa(P(2a+3,-2b-1,2c+1),1)$ to be
supported in dimension one. The inductive step on $c$ works analogously.
\end{proof}

We now turn to the case where $b\leq \min(a,c)$.

\begin{lemma}
\label{lemma:RkLowerBound}
Let $K$ be the pretzel knot $K=P(2a+1,-2b-1,2c+1)$ with $a,b,c\geq 0$
and $b\leq \min(a,c)$.
Then, $$\rk~\HFKa_{\leq 1}(K,1)\leq b(b+1).$$
\end{lemma}
 
\begin{proof}
This is proven by induction on $a$ and $c$, starting with the basic
case where $\min(a,c)=b$. Suppose for concreteness that $b=c$. Then,
it is easy to see that $\rk \HFKa(P(2a+1,-2b-1,2b+1),1)$ is
independent of $a$: the middle term in the skein exact sequence
(Equation~\eqref{eq:SkeinPretz}) vanishes: it corresponds to $\HFKa$
of the two-component unlink, whose knot Floer homology is supported in
filtration level zero. Thus, it can be calculated in the case where
$a<b$, so Proposition~\ref{prop:EasyPretzel} applies, proving that
$\HFKa(K,1)$ is supported entirely in dimension one, where its rank is
precisely $b(b+1)$.

For the inductive step, suppose we are increasing $a$ by one, and
apply Equation~\eqref{eq:SkeinPretz}. Now, the middle term
$T_{2,2(c-b)}$ has signature $+1$, and hence $\HFKa(T_{2,2(c-b)},1)$
is supported in dimension $3/2$. In particular, the map $$H\colon
\HFKa_{\leq 1}(P(2a+1,-2b-1,2c+1),1)\longrightarrow \HFKa_{\leq 1}
(P(2a+3,-2b-1,2c+1),1)$$ is surjective, providing the inductive step
for increasing $a$. Increasing $c$ follows similarly.
\end{proof}

To proceed, we use a decorated knot projection for
$P(2a+1,-2b-1,2c+1)$, and consider the multi-filtration on states from
Section~\ref{sec:Heegaards}.

Specifically, choose a knot projection for $P(2a+1,-2b-1,2c+1)$ whose
corresponding black graph consists of two vertices with degree three,
connected by three linear strands of edges $\{x_i\}_{i=1}^{2a+1}$,
$\{y_j\}_{j=1}^{2b+1}$, and $\{z_k\}_{k=1}^{2c+1}$. The vertex meeting
the edges $x_{2a+1}$, $y_{2b+1}$, and $z_{2c+1}$ is the distinguished
black vertex. The corresponding distinguished black region is a
triangle with three edges $P$, $Q$, and $R$, with $P$ facing the
vertex corresponding to $x_{2a+1}$, $Q$ facing from $y_{2b+1}$, and
$R$ facing $z_{2c+1}$. Our distinguished edge for the decorated knot
projection is $Q$. 

Let $A_{i,j}$ resp. $B_{i,j}$, resp. $C_{i,j}$
be the Kauffman state corresponding to the tree which
is obtained by deleting $y_i$ and $z_j$ 
resp. $x_i$ and $z_j$ resp. $x_i$ and $y_j$ from the black graph.

\begin{figure}
\mbox{\vbox{\epsfbox{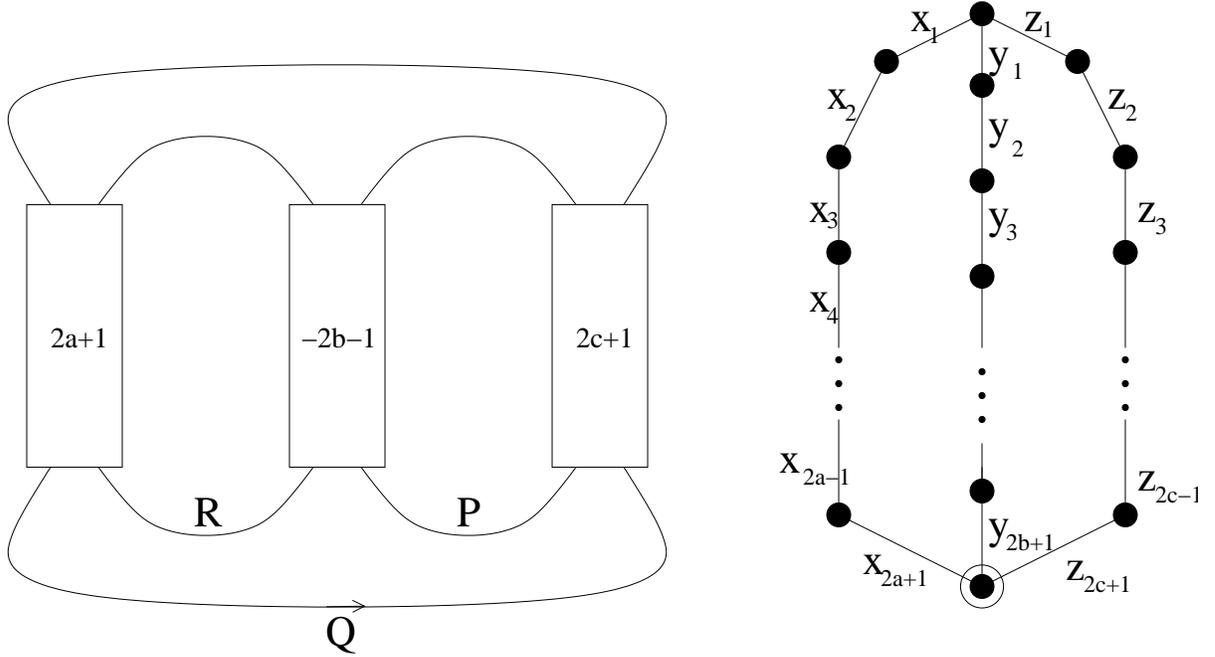}}}
\caption{\label{fig:PretzelDiag}
{\bf{Labelings for the pretzel knot P(2a+1,-2b-1,2c+1).}}
At left, we have a projection for the pretzel knot $P(2a+1,-2b-1,2c+1)$
where, of course, the labeled rectangles represent tangles with a specified
number of twists in them. The edges $P$, $Q$, $R$ are labeled here.
At the right, we have the corresponding black graph (for one of the two
colorings). The circled vertex corresponds to the black region distinguished
by the edge $Q$.}
\end{figure}

\begin{lemma}
\label{lemma:DimGenerators}
For the pretzel knot $K=P(2a+1,-2b-1,2c+1)$ with $a,b,c\geq 0$, 
the generators of $\CFKa(K,1)$ in dimension one are
$$\{A_{2i,2j+1}\}_{\begin{tiny}
			\begin{array}{l}
			1\leq i\leq b \\
			1\leq j\leq c
		\end{array}\end{tiny}}, 
\{C_{2i,2j+1}\}_{\begin{tiny}
			\begin{array}{l}
			1\leq i\leq a \\
			1\leq j\leq b
			\end{array}\end{tiny}}$$
and the generators in dimension two are of the form
$$\{B_{2i,2j+1}\}_{\begin{tiny}
			\begin{array}{l}
				1\leq i\leq a \\
				1\leq j\leq c
			\end{array}
		   \end{tiny}}.$$
There are no other generators of $\CFKa(K,1)$.
\end{lemma}

\begin{proof}
This follows at once from the above diagram. 
\end{proof}

\begin{lemma}
\label{lemma:MFPretzel}
For the marked edges $P$ and $R$, we have that
\begin{eqnarray*}
M_{B_{2i, 2j+1}}(P)-M_{A_{2s,2t+1}}(P)
&=& 2(j-s-t,j-s-t), \\
M_{B_{2i, 2j+1}}(R)-M_{C_{2s,2t+1}}(R)
&=&2(i-s-t,i-s-t).
\end{eqnarray*}
\end{lemma}

\begin{proof}
This, too, follows quickly from the diagram.
\end{proof}

\begin{lemma}
\label{lemma:RkUpperBound}
For the pretzel knot $P(2a+1,-2b-1,2c+1)$ with $a,b,c\geq 0$, we have that
$\rk \HFKa_1(P(2a+1,-2b-1,2c+1),1)\geq b\cm (b+1)$.
\end{lemma}

\begin{proof}
According to Lemma~\ref{lemma:MFPretzel} (together with
Proposition~\ref{prop:Domains}), the generators $A_{2s,2t+1}$ (where
$0\leq s \leq b$ and $0\leq t \leq c$) with $s+t>2c+1$, of which
there are $b(b+1)/2$, all lie in the cokernel of the boundary
operator.  Similarly, the generators $C_{2s,2t+1}$ 
(where $0\leq s\leq a$ and $0\leq t \leq b$) with $s+t>2a$, of which 
there are another $b(b+1)/2$, all lie in the cokernel of the boundary operator.
Since there are no generators in dimension zero (according
to Lemma~\ref{lemma:DimGenerators}), the stated bound follows.
\end{proof}

\vskip.2cm
\noindent{\bf{Proof of Theorem~\ref{thm:PretzCalc}.}}
The case where $b\geq \min(a,c)$ is established in
Proposition~\ref{prop:EasyPretzel}.  Together, Lemmas~\ref{lemma:RkLowerBound} and
~\ref{lemma:RkUpperBound} show that 
$\rk\HFKa(P(2a+1,-2b-1,2c+1))=b(b+1)$. The rest of the theorem now
follows at once from Lemma~\ref{lemma:DimGenerators}.
\qed
\vskip.2cm

To prove Corollary~\ref{cor:NoSeifSurgPretz}, we need the following result,
which closely follows~\cite{SurgSeif}:

\begin{prop}
\label{prop:NoSeifSurg}
Let $K$ be a knot with $\deg\HFKa(S^3,K)=1$. Then, if 
\begin{eqnarray*}
\rk\HFKa_{\ev}(S^3,K,1)\geq 2 &{\text{and}}&
\rk\HFKa_{\odd}(S^3,K,1)\geq 1,
\end{eqnarray*} then 
no integral surgery of $S^3$ along $K$ is a Seifert fibered space.
\end{prop}
 
\begin{proof}
By reflecting the knot if necessary, we can assume that $S^3_p(K)$ is
Seifert fibered for some $p\geq 0$. 
 
According to 
Lemma~\ref{SurgSeif:lemma:RelateLemma} of~\cite{SurgSeif}, 
$$\HFpRedEv(S^3_0(K),0)\cong \HFKa_{\odd}(S^3,K,1).$$
This completes the case where $p=0$, see for example Theorem~\ref{SurgSeif:thm:HFpObstruction} of~\cite{SurgSeif}.

As in the proof of
%\smargin{Check this}
Lemma~\ref{SurgSeif:lemma:RelateLemma} of~\cite{SurgSeif},
Section~\ref{Knots:sec:Relationship} of~\cite{Knots} gives 
a $\Z[U]$-submodule of $\HFp(S^3_n(K),[0])$ (for sufficiently large $n$)
which is isomorphic
to $\HFKa(S^3,K,1)$; indeed, we have a short  exact sequence:
$$
\begin{CD}
0@>>>\HFKa(S^3,K,1) @>>> \HFp(S^3_n(K),[0]) @>>> \HFp(S^3)@>>>0.
\end{CD}
$$
The above is a map of $U$-modules, and the $U$ action on $\HFKa(S^3,K,1)$
is trivial. It follows now that 
\begin{eqnarray*}
\rk \HFpRedEv(S^3_n(K))&\geq& \rk \HFKa_{\ev}(S^3,K,1)-1\\
\rk \HFpRedOdd(S^3_n(K))&\geq& \rk \HFKa_{\odd}(S^3,K,1).
\end{eqnarray*}
Considering the integer surgeries long exact sequence, it follows
that for all $p>0$
\begin{eqnarray*}
\rk \HFpRedOdd(S^3_n(K))&=& \rk \HFpRedEv(S^3_p(K)) \\
\rk \HFpRedEv(S^3_n(K))&=& \rk \HFpRedOdd(S^3_p(K)).
\end{eqnarray*}
In view of our hypotheses, then $\HFpRed(S^3_p(K))$ is non-trivial in
both even and odd degrees. On the other hand, results
from~\cite{SomePlumbs} show that for a Seifert fibered space
with $b_1(Y)=0$, $\HFpRed(Y)$ is supported in either
even or odd degrees (this is proved in
Corollary~\ref{Seifert:cor:EvenDegrees} of~\cite{SomePlumbs}).
\end{proof}

\vskip.2cm
\noindent{\bf{Proof of Corollary~\ref{cor:NoSeifSurgPretz}.}}
This is a direct consequence of Theorem~\ref{thm:PretzCalc}
and Proposition~\ref{prop:NoSeifSurg}.
\qed

\section{Knots with few crossings}
\label{sec:SmallKnots}

We give here another application of the results of Proposition~\ref{prop:Simplify},
showing that the Floer homology groups of all but two of the knots
with nine or fewer crossings behave like the Floer homology of
alternating knots. The two counterexamples to this are the
$(3,4)$-torus knot (which appears in the tables under the name
$8_{19}$) and a certain nine-crossing knot $9_{42}$. In fact, the knot
Floer homologies of these two knots have been determined in
Theorem~\ref{NoteLens:thm:FloerHomology} of~\cite{NoteLens}
and Proposition~\ref{Knots:prop:NineFortyTwo}
of~\cite{Knots} respectively, where it is shown that:
\begin{eqnarray*}
\HFKa(8_{19},i)&\cong& \left\{\begin{array}{ll}
\Z_{(0)} & \text{if $i=3$} \\
\Z_{(-1)} & \text{if $i=2$} \\
\Z_{(-4)}\hskip1.1cm & \text{if $i=0$} \\
\Z_{(-5)} & \text{if $i=-2$} \\
\Z_{(-6)} & \text{if $i=-3$} \\
0 & \text{otherwise,}
\end{array}\right. \\ \\
\HFKa(9_{42},i)&\cong&
\left\{\begin{array}{ll}
\Z_{(1)} & {\text{if $i=2$}}  \\ 
\Z^2_{(0)} & {\text{if $i=1$}} \\
\Z^2_{(-1)} \oplus \Z_{(0)} & {\text{if $i=0$}} \\ 
\Z^2_{(-2)} & {\text{if $i=-1$}} \\ 
\Z_{(-3)} & {\text{if $i=-2$}} \\ 
0 & {\text{otherwise.}}  
\end{array}\right.
\end{eqnarray*}
(Note that the standard knot tables do not distinguish a knot from its mirror.
For the above statements, we have chosen the versions of the knots whose
signature is negative.)

\begin{theorem}
\label{thm:SmallKnots}
Except for the knots $8_{19}$ and $9_{42}$, 
any other knot $K$ admitting a projection with nine or fewer crossings has the property that
\begin{equation}
\label{eq:CalcSmall}
\HFKa(K,i)\cong \Z^{|a_i|}_{(i+\frac{\sigma}{2})},
\end{equation}
where here $\sigma$ denotes the signature of the knot $K$, and the $a_i$ are the coefficients
of its symmetrized Alexander polynomial.
\end{theorem}

\begin{proof}
Of course,
for alternating knots, the theorem follows from~\cite{AltKnots}. 
Now, there there are only nine non-alternating knots to consider here
according to standard knot tables, see for example~\cite{BurdeZieschang}.
One of these, the knot $9_{46}$, which is the pretzel knot $P(-3,3,3)$, will
be handled separately. We illustrate 
distinguished edges for knot projections for the remaining
eight knots in 
Figure~\ref{fig:Knots} (but dropping orientations).

\begin{figure}
\mbox{\vbox{\epsfbox{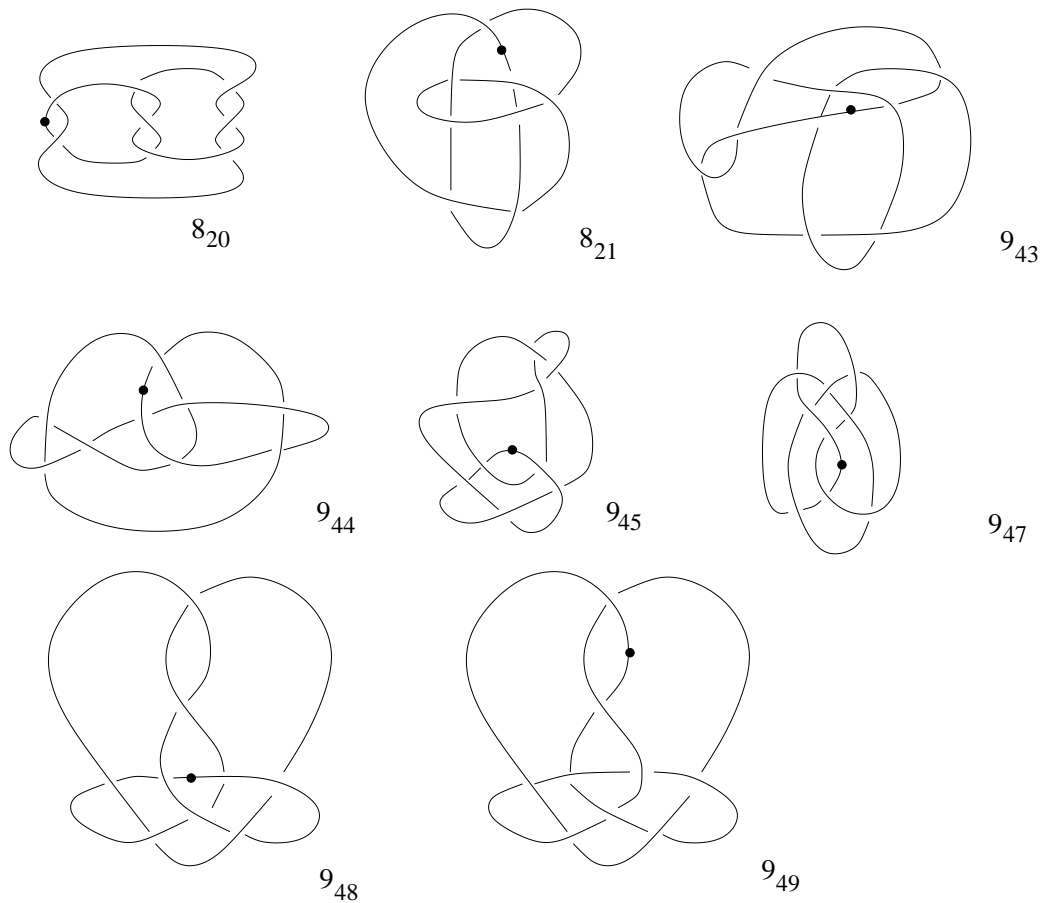}}}
\caption{\label{fig:Knots}
{\bf{Decorated knot projections for small knots.}} 
We have illustrated here knots with nine or fewer crossing which do not admit alternating
projections, except for $8_{19}$ (the $(3,4)$ torus knot), $9_{42}$, and $9_{46}$ (the pretzel
knot $P(-3,3,3)$).}
\end{figure}

Now, of these eight, we consider $9_{43}$ separately as well. For the
remaining seven knots, it is straightforward to see that in each
filtration level, all of the essential states have the same absolute
grading. Indeed, calculating these absolute gradings, one can readily
verify that for these knots, the essential states with filtration
level $i$ all have absolute grading $i+\sigma/2$. In view of
Proposition~\ref{prop:Simplify}, the theorem then follows for these
seven knots.

For the case of $9_{43}$, a direct analysis using the illustrated
decorated knot projection verifies Equation~\eqref{eq:CalcSmall} for
all $i<0$, and hence also for all $i\neq 0$, in view of the symmetry
of $\HFKa$, Equation~\eqref{eq:Symmetry}.  In the case where $i=0$,
now, we claim that there are three generators, one in dimension $-1$,
and two in dimension $-2$. In fact, a closer look at the states
reveals that the essential state $X$ in dimension $-1$ can be
connected to an essential state $Y$ in dimension $-2$ by a homotopy
class $\phi$ whose associated domain is an octagon (with multiplicity
$+1$, missing the reference point $z$). Compare  with
the genus four Heegaard diagram of $S^3$ pictured in
Figure~\ref{fig:Octagons}, where there are three generators for $S^3$,
$A_1, A_2, B$, with homotopy classes connecting $A_1$ resp. $A_2$ to
$B$ represented by octagons. It now follows easily that for a homotopy
class $\phi$ whose domain is an octagon,
$$\#\left(\ModFlow(\phi)/\R\right) = \pm 1.$$

Hence, we have that the boundary operator in $\CFKa(9_{43},0)$
is non-trivial, and indeed
that the homology in filtration level $0$ is given by $\Z_{(-2)}$,
completing the verification of Equation~\eqref{eq:CalcSmall} for
$9_{43}$.

\begin{figure}
\mbox{\vbox{\epsfbox{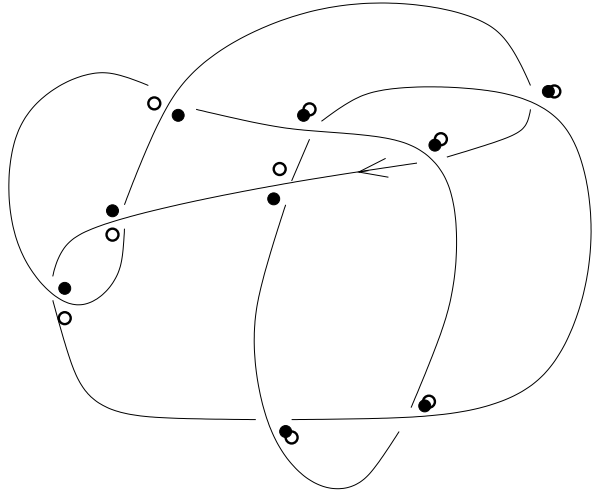}}}
\caption{\label{fig:9s43}
{\bf{A differential for $9_{43}$.}} 
We have illustrated here two of the essential states $X$ and $Y$ for the indicated decorated
knot projection (where the distinguished edge 
is the one containing the arrow). The 
state $X$ is represented by the collection of dark circles, while $Y$ is represented
by the hollow circles. Moreover, $X$ and $Y$ are in dimensions $-1$ and $-2$ respectively, 
and it is easy to see that the domain of the homotopy class connecting $X$ to $Y$ is 
an octagon.}
\end{figure}

\begin{figure}
\mbox{\vbox{\epsfbox{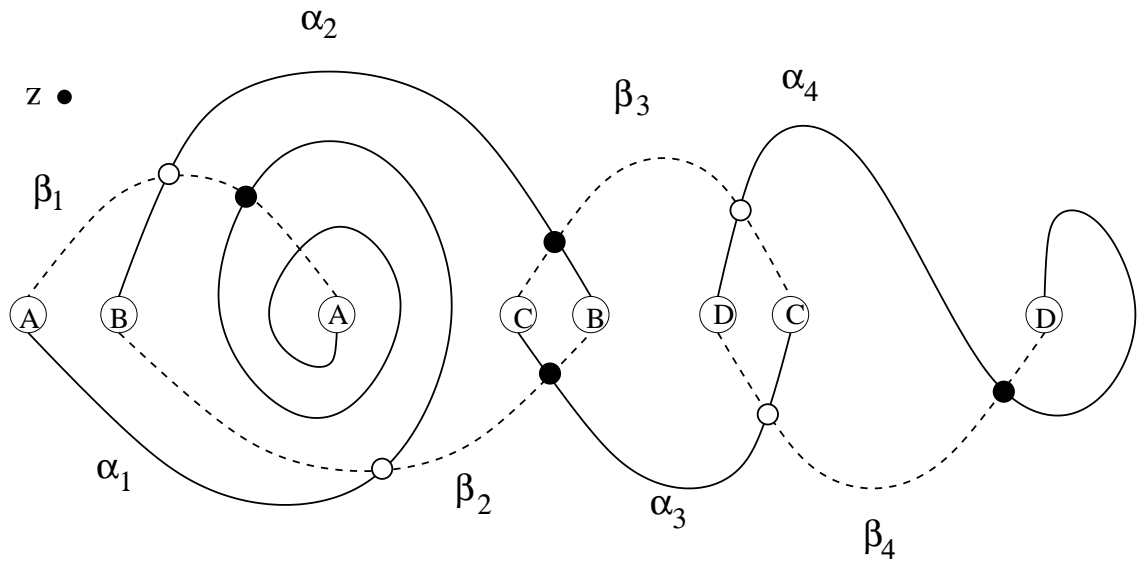}}}
\caption{\label{fig:Octagons}
{\bf{Octagons.}}  In this genus four Heegaard diagram for $S^3$, there
are three generators for $\CFa(S^3)$. Two of them are indicated here
-- one by the unmarked solid circles (call it $X$), the other by the
unmarked hollow circles (call that $Y$).  It is easy to find an
octagonal domain $\cald(\phi)$ with $n_z(\phi)=0$ which connects $X$ to $Y$ (and indeed
there is another octagonal domain connecting the other intersection
point $X'$ to $Y$). This forces $\#\ModFlow(\phi)/\R=\pm 1$.}
\end{figure}

Finally, we turn to the pretzel knot $P(-3,3,3)$.
As in Section~\ref{sec:PretzCalc}, we fit this into a skein exact sequence
$$
\begin{CD}
...@>>> \HFKa(P(-3,3,3)) @>{F}>> \HFKa(U_2) @>{G}>> \HFKa(P(-3,3,1))@>{H}>>...
\end{CD}
$$
where here $U_2$ is the unlink with two components (this is $T_{2,0}$ in the 
notation from Section~\ref{sec:PretzCalc}). By using the action of the
homology class which links, it is easy to see that $\HFKa(P(-3,3,3))\cong \HFKa(P(-3,3,1))$.
Note that $P(-3,3,1)=P(-1,-1,5)=P(1,1,3)$, which alternates. 
\end{proof}

\commentable{
\bibliographystyle{plain}
\bibliography{biblio}
}

\end{document}